\documentclass[a4paper,12pt,draft]{article}
\usepackage{amsmath,amsfonts}
\usepackage{color}
\definecolor{Blue}{rgb}{0.1,0.1,0.9}
\definecolor{Red}{rgb}{0.9,0.1,0.1}
\definecolor{green}{rgb}{0,1,0}

\voffset - 2 cm
\hoffset - 1 cm
\newcommand{\cl}[1]{\overline{#1}}
\newcommand{\expr}[1]{\left( #1 \right)}
\newcommand{\set}[1]{\left\{#1 \right\}}
\newcommand{\ex}{\mathbf{E}}
\newcommand{\pr}{\mathbf{P}}
\newcommand{\R}{\mathbf{R}}
\newcommand{\Rd}{{\R^d}}
\newcommand{\Rdwa}{{\R^2}}
\newcommand{\ind}{\mathbf{1}}
\newcommand{\wt}{\widetilde}
\newcommand{\pim}{P^*}

\newcommand{\proof}{\noindent \textsc{Proof: }}
\newcommand{\proofof}[1]{\noindent \textsc{Proof of #1: }}
\newcommand{\qed}{\ensuremath{\hspace{10pt}\square}}

\newtheorem{theorem}{Theorem}
\newtheorem{proposition}{Proposition}
\newtheorem{lemma}{Lemma}
\newtheorem{corollary}{Corollary}
\newtheorem{definition}{Definition}
\newtheorem{remark}{Remark}

\DeclareMathOperator{\diam}{diam}
\DeclareMathOperator{\dist}{dist}
\DeclareMathOperator{\ro}{RO}
\DeclareMathOperator{\supp}{supp}

\begin{document}
\sloppy


\title{Estimates and structure of $\alpha$-harmonic functions}
\author{Krzysztof Bogdan\thanks{Supported by KBN 
grant 1 P03A 026 29
and RTN 
contract 
HPRN-CT-2001-00273-HARP}, Tadeusz Kulczycki\thanks{Supported by KBN 
grant 1 P03A 020 28
and RTN 
contract 
HPRN-CT-2001-00273-HARP}, Mateusz Kwa\'snicki\thanks{Supported by KBN 
grant 1 P03A 020 28 and RTN 
contract 
HPRN-CT-2001-00273-HARP }}
\date{3/19/2007}
\maketitle


\begin{abstract}
We prove a uniform boundary Harnack inequality for nonnegative
harmonic functions of the fractional Laplacian on arbitrary open set
$D$.
This yields a unique representation of such functions as integrals 
against measures on $D^c\cup \{\infty\}$ satisfying an integrability condition.
The corresponding Martin boundary of $D$ is a subset of the Euclidean
boundary determined by an integral test.
\end{abstract}

\footnotetext{2000 {\it MS Classification}\/:
31C35, 60J50 (Primary), 31B05, 60G51 (Secondary).\\
{\it Key words and phrases}\/: boundary Harnack inequality, Martin representation, stable process.}


\section{Main results and introduction} \label{sec:introduction}

Let $d = 1, 2, \ldots$, and $0 < \alpha < 2$.  
The boundary Harnack principle (BHP) for nonnegative harmonic
functions of the fractional Laplacian {on $\Rd$}
\begin{equation} \label{eq:laplacian}
  \Delta^{\alpha / 2} \varphi(x) =
 \lim_{\varepsilon \rightarrow 0^+} 
    \int\limits_{|y-x| > \varepsilon}
    [\varphi(y) - \varphi(x)]\,\nu(x,y) dy\,,
\end{equation}
was proved for Lipschitz domains in 1997 in \cite{bib:b:bhp} (compare Theorem~\ref{th:bhp} below).  
Here 
$$\nu(x,y)={\cal A}_{d,-\alpha}|y-x|^{-d-\alpha}\,,$$ 
${\cal A}_{d,\gamma}=\Gamma((d-\gamma)/2)/(2^{\gamma}\pi^{d/2}|\Gamma(\gamma/2)|)$
for $-2<\gamma<2$, and, say, $\varphi \in C^\infty_c(\Rd)$.
BHP was extended to all open sets in 1999 in \cite{bib:sw}, 
with the constant in the estimate depending on local geometry of their
boundary.
The question whether the constant may be chosen independently 
of the domain, or {\it uniformly}, was since open.

In what follows $D$ is an arbitrary nonempty open subset of
$\Rd$ (a domain). 
Let $G_D$ be the Green function of $D$ for $\Delta^{\alpha / 2}$
(\cite{bib:La}, \cite{bib:bh}, \cite{bib:po}). 
We define the Poisson kernel of $D$:
\begin{equation} \label{eq:poisson:definition}
  P_D(x, y) =
  \int_D G_D(x, v) \nu(v, y) \, dv
  \, , \quad x \in \Rd \, , \; y \in D^c \,.
\end{equation}
By a calculation of M. Riesz (see \cite{BGR}, \cite{bib:Rm}), for the ball $B_r = \{x\in \Rd\,:\;
|x|<r\}$ we have
\begin{equation} \label{eq:poisson:ball}
  P_{B_r}(x, y) =
  {\cal C}_{d, \alpha} \expr{\frac{r^2 - |x|^2}{|y|^2 - r^2}}^{\alpha / 2}
    \frac{1}{|x - y|^d}
  \, , \quad x \in B_r\, , \; y \in B_r^c \, ,
\end{equation}
where ${\cal C}_{d,\alpha}=\Gamma(d/2)\pi^{-1-d/2}\sin (\pi \alpha/2)$.
Note that if $x$ and $y$ are not too close then
$P_{B_r}(x, y) \approx (r^2-|x|^2)^{\alpha/2}
\cdot (|y|^2-r^2)^{-\alpha/2}|y|^{-d}$ at $\partial B_r$.
Similar approximate factorization of
general $P_D$ underlies the following result
which is equivalent to the uniform BHP (UBHP) for $\Delta^{\alpha/2}$
(see also Theorem~\ref{th:bhp} and Remark~\ref{rem:gh} below).
%
\begin{theorem} \label{th:bhp:ball} {\rm (UBHP)}
There is a constant $C_{d, \alpha}$, depending only on $d$ and
$\alpha$, such that 
\begin{equation} \label{eq:bhp}
P_D(x_1, y_1)P_D(x_2, y_2) \leq
C_{d, \alpha}\, P_D(x_1, y_2)P_D(x_2, y_1)\, ,
\end{equation}
whenever $r>0$, $x_1, x_2 \in D \cap B_{r/2}$ and
$y_1, y_2\in D^c \cap B_r^c$.
\end{theorem}
Consider the following auxiliary function
\begin{equation}
  \label{eq:ds}
s_D(x) = \int_\Rd G_D(x, v) dv \,.
\end{equation}
We will say that $y\in \Rd$ 
is {\it accessible}\/ from $D$ when
\begin{equation}\label{eq:accessibility}
  \int_\Rd s_{D \cap B(y, 1)}(v) \, \nu(v, y) dv = \infty \,,
\end{equation} 
or
{\it inaccessible}\/ when 
\begin{equation} \label{eq:inaccessibility}
  \int_\Rd s_{D \cap B(y, 1)}(v) \, \nu(v, y) dv < \infty \,.
\end{equation}
The {\it point at infinity}\/ is called {\it accessible}\/ for $D$ if $s_D(x)=\infty$ for
all $x\in D$, and it is called {\it inaccessible}\/ otherwise. Accessibility
of a given point from $D$
means that $D$ is rather large near the point, see (\ref{eq:mons}),
(\ref{eq:ginaccessible}) and the discussion at the end of the paper.

We consider the set $\partial_* D$ of {\it limit}\/ points of $D$:
we let $\partial_* D=\partial D$  if $D$ is bounded 
and $\partial_* D=\partial D\cup \{\infty \}$ if $D$ is unbounded.
For unbounded $D$, $D\ni v\rightarrow \infty$ means that
$v\in D$ and $|v|\rightarrow \infty$.
We also let $D^* = D \cup \partial_* D$.

Theorem~\ref{th:bhp:ball} 
and Theorem~\ref{th:bhp} yield estimates of
$G_D$, $P_D$ and $s_D$
at $\partial_* D$.
Furthermore, the uniformity of the constant 
yields results on the {\it limits}\/ of ratios of these kernel functions.
Our main result in this direction is Lemma~\ref{th:oscillation} below.
An important consequence of the result is the existence of the Martin kernel. 
For {\it Greenian}\/ {$D \subset \Rd$} we
fix an arbitrary reference point $x_0\in D$ and 
we define the Martin kernel of $D$:
\begin{equation} \label{eq:martin:definition}
  M_D(x, y) =
  \lim_{D \ni v \rightarrow y} \frac{G_D(x, v)}{G_D(x_0, v)}
  \, , \quad x \in \Rd \,,\;  y\in \partial_* D\,.
\end{equation}
%
\begin{theorem} \label{th:martin}
The limit in (\ref{eq:martin:definition}) exists. $M_D(x, y)$ is
$\alpha$-harmonic in $x$ on $D$ {with zero outer charge
  on $D^c$} if and
only if $y$ is accessible from $D$. 
If $y\in \partial D$ is inaccessible then
$M_D(x, y) = P_D(x, y)/ P(x_0,y)$. 
If $\infty\in \partial_*D$ is inaccessible then
$M_D(x,\infty)=s_D(x)/s_D(x_0)$.
\end{theorem}
We define
$\partial_M D=\{y\in \partial_* D\,:\, y \;\mbox{is accessible from } D\}$ and 
$D^c_M=\{y\in D^c\,:\, y \; \mbox{is inaccessible from } D \}$.
The kernels $M_D(\cdot,y)$, $y\in \partial_M D$, and
$P_D(\cdot,y)$, $y\in D_M^c$,  may be used to describe the structure of nonnegative functions 
harmonic for $\Delta^{\alpha / 2}$ on $D$, or $\alpha$-harmonic
(a detailed discussion of the notion of $\alpha$-harmonicity is given
in Section~\ref{sec:harmonicity}).
\begin{theorem} \label{th:representation}
Let $D$ be Greenian. For every function $f \geq 0$ on $D$ which is $\alpha$-harmonic
in $D$ with outer charge $\lambda\geq 0$ on $D^c$ there
is a unique  measure
$\mu\geq 0$ on $\partial_M D$, such that
\begin{equation} \label{eq:representation}
  f(x) =
\int_{D^c}
P_D(x, y) \lambda(dy) + \int_{\partial_M D} M_D(x, y) \mu(dy)
  \, , \quad x \in D \,.
\end{equation}
\end{theorem}
As a part of the statement we have that 
$|\mu| < \infty$,
\begin{equation} \label{eq:lambda}
  \int_{D^c} P_D(x_0, y) \lambda(dy) < \infty\,,
\end{equation}
and $\lambda(\partial_M D)=0$, so that $\lambda$ is concentrated on
$D^c_M$, see (\ref{eq:obM}).
We remark that for non-Greenian $D$ every such $f$ is constant on $D$,
see Lemma~\ref{l:hc}. For a stronger statement of the
uniqueness, see Remark~\ref{rem:u} below.

The first integral in (\ref{eq:representation}) reflects the fact
that $\Delta^{\alpha/2}$ is a 
{\it nonlocal}\/ integro-differential operator, allowing
  for a direct integral-type influence between distant points $x$ and $y$ in the domain of a
  function, see (\ref{eq:laplacian}). 
In particular the role of the boundary condition in the 
Dirichlet problem of the classical potential theory is now played by 
a measure (``outer charge'') supported on the complement of the domain.
Here the generic example is the Poisson kernel $P_D(x,y)$ equipped
with the Dirac measure at $y\in D^c$.
We remark that the restriction to harmonicity of only genuine functions
would seriously handicap the theory because 
the limit of a locally bounded pointwise convergent
sequence of genuine functions which are harmonic for
$\Delta^{\alpha/2}$ on a given domain may fail to be
a genuine $\alpha$-harmonic function itself 
(see the concluding remark
in \cite{bib:ms},
related to inaccessible boundary points of $D$).

Theorem~\ref{th:martin} and Theorem~\ref{th:representation} contrast
sharply with the corresponding results in the classical potential
theory (\cite{bib:AG}, \cite{bib:pi}), because they are
  more explicit, and also because the {\it classical}\/ Martin kernel is
  always harmonic, which is no longer the case here.
 We refer the interested reader to \cite{bib:KW} and
  \cite{Br} for a general account on Martin compactification and
  representation. We also refer to the paper \cite{bib:hw},
  which identifies the classical Martin boundary of Lipschitz domains 
with their Euclidean boundary
(see  \cite{bib:pi} or \cite{bib:Ba} for further references). 
We see that the fractional Laplacian enjoys a similar description 
  in an arbitrary domain. 

The role of BHP in explicit determination of the Martin boundary in
the classical potential theory is well recognized, see
recent \cite{bib:A} and \cite{bib:A2} (see also \cite{bib:BB}, 
\cite{bib:Ba}, and \cite{bib:AG} for more references).
The present straightforward derivation of the Martin representation 
is modeled after \cite{bib:b:repr}.
The role of BHP in estimating the Green function and
studying Schr\"odinger-type operators is also well
understood. For more information on typical applications 
we refer the reader to \cite{Bjmaa2000}, \cite{bib:TJe}, \cite{bib:BBC},
\cite{bib:CK}, \cite{BB1}, \cite{BBpms2000}, \cite{CS:FK}, see
also \cite{bib:Hg} for a general perspective.

Our theorems complete and extend in several directions
part of the results of \cite{bib:b:bhp}, \cite{bib:ki}, \cite{bib:b:repr}, \cite{bib:cs}, 
\cite{bib:sw}, \cite{bib:ms}. In
particular, Theorem~\ref{th:representation} was first proved for Lipschitz domains in
\cite{bib:b:repr} and \cite{bib:cs}, and for $\kappa$-fat domains in
\cite{bib:sw}. For these domains all the boundary points are accessible, which
influenced the methods of these papers. The first example of what we
coin inaccessible 
boundary point was given in \cite{bib:ms}. 
Our main technical results, Lemma~\ref{th:factorization} and
Lemma~\ref{th:oscillation},
develop the ideas of \cite{bib:b:bhp} (see also the references in
\cite{bib:b:bhp}) and \cite{bib:sw}.

The paper is primarily addressed to the readers interested in the
potential theory of nonlocal operators. The theory presently undergoes a rapid
development, see \cite{bib:J} and the references given there. The outline and notions which we
propose below may likely apply to kernel functions of such operators
and the corresponding nonnegative
harmonic functions quite generally, except for our treatment of the
point at infinity, which is based on Kelvin
transformation and therefore is very specific to the present context.
Technically, the development hinges on
Lemma~\ref{th:factorization} and Lemma~\ref{th:oscillation} below, and
extensions of these should be sought for in the more general
settings. A certain role is also played by (\ref{eq:omega:density}).

Our development is based on M. Riesz'
formulas (\ref{eq:poisson:ball}) and (\ref{wfg}) for the Poisson
kernel and the Green function of the ball, and
general properties of the Green function and
harmonic measure of arbitrary domains, most notably (\ref{eq:green:twosets}). Here our
references are \cite{bib:po}, \cite{bib:La}, and \cite{bib:bh}. 
The reader familiar with the potential theory of Markov processes will
notice the relationship of (\ref{eq:green:twosets}) to the strong Markov property of the
isotropic $\alpha$-stable L\'evy process $\left\{X_t,\,t\geq 0\right\}$ in ${\bf R}^d$ with the
L\'evy measure $\nu(0,x)dx$, see \cite{bib:bz}, \cite{bib:Sa}. 
Indeed, probabilistic interpretations and references are our primary
source of motivation, as seen from the discussion of
such interpretations given at the end of the paper.
In the main body of the paper we strive, however, to give elementary and purely {\it analytic}\/ definitions
and proofs, with a notable exception made for the {\it probabilistic}\/ proof of Lemma~\ref{th:omega:twosets}. 

The remainder of the paper is organized as follows. 
In Section~\ref{sec:preliminaries} we give preliminary definitions and results.
In Section~\ref{sec:bhp} we prove Theorem~\ref{th:bhp:ball} and we state 
UBHP in a more traditional form as Theorem~\ref{th:bhp}.
In Section~\ref{sec:limits} we study limits of ratios of kernel functions.
In Section~\ref{sec:harmonicity} we
define $\alpha$-harmonicity.
In Section~\ref{sec:martin} we prove joint continuity of $M_D(x,y)$
and verify Theorem~\ref{th:martin}. 
In Section~\ref{sec:harmonic} we obtain the Martin representation
(\ref{eq:representation}) along with its converse.
In Section~\ref{sec:miscelanea} we prove absolute continuity of harmonic
measure on $D^c_M$, discuss probabilistic interpretations of our results
and give examples of accessible and inaccessible boundary points. 
For instance $0$ is inaccessible for $D=\{(x,y)\in
\Rdwa:\;y>|x|^\gamma\}$ if and only if $\gamma<1$.




\section{Preliminaries} \label{sec:preliminaries}

For $x \in \Rd$ and $r > 0$ we let $|x| = \sqrt{\sum_{i=1}^d x_i^2}$,
$B(x,r) = \{y\in\R :\: |y-x|<r \}$, $B_r=B(0, r)$, and $B=B_1$. 
All the sets, functions and measures considered in the sequel 
will be Borel. 
For $U\subset \Rd$ we write $U^c=\Rd\setminus U$. If $k>0$ then
$kU=\{kx:\; x\in U\}$. 
For a measure $\lambda$ on $\Rd$, $|\lambda|$ denotes its total mass. For a
function $f$ we let $\lambda(f) = \int f d\lambda$ if the integral makes sense.
The probability measure concentrated at $x$ 
will be denoted by $\varepsilon_x$. 
For nonnegative $f$ and $g$ and a positive number $C$ we write
$f\asymp C\, g$ if $C^{-1}f\leq g\leq C f$.
{The notation $C_{a,b,\ldots,z}$ means that such constant depends
{\it only}\/ on $a,b,\ldots,z$.}
In what follows $U$ will be an arbitrary domain.
We will say that $U$ is Greenian if $G_U(x,v)$ 
is finite almost everywhere on $U\times U$. 
$U$ is always Greenian when $\alpha<d$. 
If $\alpha\geq d=1$, then $U$ is Greenian if and only if $U^c$ is
non-polar. In particular, if $\alpha>d=1$, then $U$ is Greenian unless
$U={\bf R}$. 
The Green function and the harmonic measure 
of the fractional Laplacian are defined in 
\cite[Theorem IV.4.16, pp. 229, 240]{bib:La}, see also \cite{BGR},
\cite[pp. 191, 250, 384]{bib:bh},
\cite{bib:Kg},
and \cite{bib:po}, \cite{bib:bz} for the case of dimension one.
We will briefly indicate the following crucial properties.
If $U$ is Greenian then 
\begin{equation} \label{eq:green:definition}
  \int_\Rd G_U(x, v) \Delta^{\alpha / 2} \varphi(v) dv =
  -\varphi(x)
  \, , \quad x \in \Rd \, , \, \varphi \in C_c^\infty(U) \,.
\end{equation}
Furthermore, $G_U(x, v) = G_U(v, x)$ for $x, v \in \Rd$
(\cite[p. 285]{bib:La}).
For example, if $\alpha<d$ then
  (\ref{eq:green:definition}) is satisfied for $U=\Rd$ by the Riesz kernel: 
\begin{equation}\label{eq:potr}
G_\Rd(x,y)=
{\mathcal A}_{d,\alpha}|y-x|^{\alpha-d}\,,\quad x,y\in \Rd\,,
\end{equation}
see \cite[(1.1.12')]{bib:La}, and the harmonic measure, $\omega^x_U$, is defined as the unique
(\cite[p. 245]{bib:La}, \cite{bib:bz})
subprobability measure (probability measure if $D$ is bounded) 
concentrated on $D^c$ such that $\int_{\Rd} G_\Rd(z,y)
\omega^x_U(dz)\leq G_\Rd(x,y)$ for all $y\in \Rd$, and 
\begin{equation}
  \label{eq:dhm}
  G_\Rd(x,y)=\int_{\Rd} G_\Rd(z,y) \omega^x_U(dz) 
\end{equation}
for $y\in U^c$ except at {\it irregular}\/ points of $\partial U$.
Recall that a point $y$ is called {\it irregular}\/ for $U$ (or {\it thin}\/ for $U^c$)
if $\omega^y_U\neq\varepsilon_y$, and
it is called {\it regular}\/ otherwise, see \cite[pp. 348, 272, 353]{bib:bh}. Note that
``regularity'' means here ``regularity for the Dirichlet problem on $U$'' \cite[p. 348]{bib:bh}.
The {probabilistic} interpretation of regularity is that the first
{\it hitting}\/ time of $U^c$ for the corresponding stochastic process starting at $x$ equals
zero almost surely, see \cite[p. 277]{bib:bh}.
Note that $y$ is regular for $U$ if and only if $G_U(x,y)=0$ for 
$x\in U$ (\cite[Proposition VII.3.1]{bib:bh}, see also \cite[pp. 251, 286]{bib:La}).
If $\alpha<d$ then the Green function is given by
\begin{equation}
  \label{eq:wnfg}
  G_U(x,y)=G_\Rd(x,y)-\int_{U^c} G_\Rd(z,y) \omega^x_U(dz)\,.
\end{equation}
For a full discussion we refer the reader to 
\cite{bib:La}, \cite{bib:bh} (see also \cite{bib:po} for $1=d\leq \alpha$).
The harmonic measure
may be used to negotiate between Green functions of
two domains:
\begin{equation} \label{eq:green:twosets}
  G_D(x, v) = G_U(x, v) + \int_\Rd G_D(w, v) \omega^x_U(dw)
  \, , \quad x, v \in \Rd \, , \, \mbox{if }U \subset D\,,
\end{equation}
compare (\ref{eq:wnfg}).
By integrating \eqref{eq:green:twosets} against the Lebesgue measure
we obtain
\begin{equation} \label{eq:superharmonic}
  s_D(x) =
  s_U(x) + \int_\Rd s_D(y) \omega^x_U(dy) \,, \quad x \in  \Rd\,,\; 
U\subset D\,.
\end{equation}
Clearly,
\begin{equation} \label{eq:mons}
  s_U\leq s_D\,,\quad \mbox{if } U\subset D\,.
\end{equation}
Recall that $\supp \omega^x_U\subset U^c$, $x\in \Rd$.
If $U \subset D$ then
\begin{equation} \label{eq:omega:harmonic}
  \omega^x_D(A) =
  \omega^x_U(A) + \int_{D \setminus U} \omega^y_D(A) \omega^x_U(dy)
  \, , \quad A \subset D^c \,,
\end{equation}
in particular $G_{U}(x, v) \le G_{D}(x, v)$ and for
$A \subset D^c$, $x, v \in U$, we have 
\begin{equation}
  \label{eq:15.5}
\omega^x_{U}(A) \le 
\omega^x_{D}(A)\,.
\end{equation}
Furthermore, if
$D_1 \subset D_2 \subset \ldots$ and $D = \bigcup D_n$, then
$G_{D_n}(x, v) \uparrow G_D(x, v)$ and $\omega^x_{D_n}(\varphi) \rightarrow
\omega^x_D(\varphi)$ whenever $x, v \in D$ and $\varphi \in C_0(\Rd)$
({\it vague}\/ convergence, {\it weak}\/ convergence for bounded $D$ \cite[(4.6.6)]{bib:La}).

Let $\varphi \in C_c^\infty(\Rd)$ and let open Greenian $D'$ contain both $D$ and
the support of $\varphi$. Using \eqref{eq:green:definition} for $D$
and $D'$, \eqref{eq:green:twosets}, and Fubini we obtain
\begin{equation} \label{eq:green:identity}
  \int_D G_D(x, v) \Delta^{\alpha/2} \varphi(v) dv =
  \int_{D^c} [\varphi(y) - \varphi(x)] \omega^x_D(dy)
  \, , \quad x \in D \, , \, \varphi \in C_c^\infty(\Rd).
\end{equation}
By considering $\varphi$ supported away from $\cl{D}$, and by (\ref{eq:laplacian})
we conclude that on $(\cl{D})^c$, 
$\omega^x_D$ is absolutely
continuous with respect to the Lebesgue measure, and has density $P_D(x, y)$ 
given by (\ref{eq:poisson:definition}).
This is the Ikeda-Watanabe formula (\cite{bib:iw}):
\begin{equation}
  \label{eq:n}
\omega^x_D(A)=\int_A P_D(x,y)dy\,,\quad \mbox{if }\; \dist(A,D)>0\,.
\end{equation}
If $D'
{\supset}D$ is a Lipschitz domain (e.g.\ a ball) then 
$\omega^x_D(\partial D') \le \omega^x_{D'}(\partial D') = 0$
(\cite{bib:b:bhp}), hence
\begin{equation} \label{eq:omega:density}
  \omega^x_D(dy) = P_D(x, y) dy \; \mbox{on } D'^c
  \quad \mbox{provided } x\in D \subset D'\; \mbox{and $D'$ is Lipschitz} \,.
\end{equation}
The Green function of the ball is known explicitly: 
\begin{equation}\label{wfg}
G_{B_r}(x,v)={\cal
B}_{d,\alpha}\,|x-v|^{\alpha-d}\int_{0}^{w}\frac{s^{\alpha/2-1}}
{(s+1)^{d/2}}\,ds\,,\quad x,v \in B_r,
\end{equation}
where 
\begin{displaymath}
w=(r^2-|x|^{2})(r^2-|v|^{2})/|x-v|^{2},
\end{displaymath}
and ${\cal B}_{d,\alpha}=
\Gamma(d/2)/(2^{\alpha}\pi^{d/2}[\Gamma(\alpha/2)]^{2})$,
see \cite{BGR}, \cite{bib:Rm}.
It is also known (\cite{BBpms2000}, \cite{bib:bk}) that
\begin{equation}
  \label{eq:sBr}
  s_{B_r}(x)=\frac{{\cal C}_{d,\alpha}}{{\cal 
  A}_{d,-\alpha}}(r^2-|x|^2)^{\alpha/2}\,,\quad |x|\leq r\,.
\end{equation}


For a nonnegative measure $\lambda$ on $\Rd$ ($D^c$) we define its {\it Poisson integral}\/ on $D$,
$$
  P_D[\lambda](x) =
  \int_{D^c} P_D(x, y) \lambda(dy)
  \, , \quad x \in D \,,
$$
compare \eqref{eq:representation}. 
To simultaneously control $P_D[\lambda]$ and $\lambda$ we 
define the {\it measure}
\begin{equation}
  \label{eq:dhl}
  \pim_D[\lambda](dx) = P_D[\lambda](x)dx + \lambda(dx)\,.
\end{equation}
Thus, $\pim_D[\lambda]$ is equal to $\lambda$ on $D^c$,
and on $D$ it is absolutely continuous with respect to the Lebesgue
measure, with $P_D[\lambda]$ as the density function. Of course,
$P_D[\pim_D[\lambda]](x)=P_D[\lambda](x)$, $x\in D$. This observation
will be strengthened in (\ref{eq:cm}) below.


If $U\subset D$ and $v\in U^c$ is such that
$G_U(x,v)=0$ for $x\in \Rd$ (in particular, if $U$ is a Lipschitz
domain and $v\in U^c$ is arbitrary), then by \eqref{eq:green:twosets} we
have
\begin{equation} \label{eq:green:harmonic}
  G_D(x, v) =
  \int G_D(w, v) \omega^x_U(dw)
  \, , \quad x \in U\,.
\end{equation}
This, and (\ref{eq:poisson:harmonic}) below may be considered a mean value property.

The following sum of integrals will be important.
Consider a nonnegative function $f$ on $D$, a nonnegative measure
$\lambda$ on $D^c$ and a nonempty open set $U\subset D$. We denote
\begin{equation} \label{eq:cm}
\Omega^D_U[f,\lambda](x)
=\int_{D\setminus U}f(y)\omega^x_U(dy)
+
\int_{D^c}P_U(x,y)\lambda(dy)
  \,,\quad x \in U \,.
\end{equation}
Informally, we may think of $\Omega^D_U[f,\lambda]$ as an integral of $f+\lambda$
against the harmonic measure $\omega_U$. The delicate point of the
definition is that 
the integration over $\partial U\cap \partial D$ 
is restricted to the part of the harmonic measure which is
absolutely continuous with respect to the Lebesgue measure, with the
density function given by the Poisson kernel. 
The convention will play a role for $U$ touching $\partial D$. 
In this connection see (\ref{eq:n}),
Proposition~\ref{th:omega:poisson} and the discussion at the end of
the paper. 
\begin{lemma} \label{th:poisson:harmonic}
If $U \subset D$ and $\lambda$ is a nonnegative measure on $D^c$, then
\begin{equation} \label{eq:poisson:harmonic}
  P_D[\lambda](x) = \Omega^D_U[P_D[\lambda],\lambda](x)\,,\quad x \in U \,.
\end{equation}
\end{lemma}
\proof
Let $x \in U$, $y \in D^c$.
By integrating \eqref{eq:green:twosets} against $\nu(v,y)dv$ on $\Rd$,
and (\ref{eq:poisson:definition}), we get
\begin{equation}
  \label{eq:harm}
P_D[\varepsilon_y](x)=P_D(x, y)=P_U(x, y)+ \int P_D(z, y) \omega^x_U(dz)
= \Omega^D_U[P_D[\varepsilon_y],\varepsilon_y](x)\,.  
\end{equation}
The case of general $\lambda\geq 0$ follows from Fubini-Tonelli theorem.
\qed

The next two lemmas are versions of Harnack inequality, see also
Remark~\ref{rem:har}.
\begin{lemma}\label{lem:Hi}
  If $\lambda\geq 0$ and $x_1, x_2\in B_{r}\subset B_s\subset D$
  then 
  \begin{equation}
    \label{eq:nH}
    P_D[\lambda](x_1)\leq \left(\frac{1+r/s}{1-r/s}\right)^{d} P_D[\lambda](x_2)\,.
  \end{equation}
\end{lemma}
\proof
By (\ref{eq:poisson:ball}) we have
$P_{B_s}(x_1,z)\leq
(1+r/s)^{d}(1-r/s)^{-d}P_{B_s}(x_2,z)$ if $|z|\geq
{s}$. Using the second equality in (\ref{eq:harm}) with $U=B_s$, (\ref{eq:omega:density}), and
(\ref{eq:poisson:ball}), we prove the result. \qed

\begin{lemma}\label{lem:Hi2}
  If $x_1,x_2\in D$ then there is $c_{x_1,x_2}$ such that
  for every $\lambda\geq 0$
  \begin{equation}\label{eq:nH2}
    P_D[\lambda](x_1)\leq c_{x_1,x_2} P_D[\lambda](x_2)\,.
  \end{equation}
\end{lemma}
\proof
If $x_1,x_2\in B_r\subset B_{2r}\subset D$ for some $r>0$ then we are done by
Lemma~\ref{lem:Hi} with $c=c_{x_1,x_2}$ depending only on $d$.
Assume that $B(x_1,2r)\subset D$, $B(x_2,2r)\subset D$, $B(x_1,2r)\cap
B(x_2,2r)=\emptyset$ for some $r>0$, and consider (\ref{eq:harm}) with
$U=B(x_1,r)$. Let $y\in D^c$. By (\ref{eq:omega:density}) and the
first part of the proof we obtain 
$P_D(x_1,y)\geq \int_{B(x_2,r)}cP_D(
{x_2},y)P_{B_r}(0,x-x_1)dx$.
\qed

\noindent
If $K\subset D$ is compact and $x_1,x_2\in K$ then $c_{x_1,x_2}$ in Harnack's
inequality above
may be so chosen to depend only on $K$, $D$, and $\alpha$, because $r$ in the above proof may be chosen
independently of $x_1,x_2$.
Note that $D$ and $K$ may be disconnected.

\begin{remark}\label{r:c}
{\rm
If $\lambda\geq 0$ and $P_D[\lambda](x)$ is finite (positive) for
some $x\in D$, then it is locally bounded from above (below, resp.) for all $x\in D$.
This follows from Lemma~\ref{lem:Hi2}. 
Note that if (\ref{eq:lambda}) holds then $P_D[\lambda]$ is finite and
locally uniformly Lipschitz continuous on $D$, a consequence of
(\ref{eq:nH}).
}  
\end{remark}
The following well-known result is given for the reader's convenience.
\begin{lemma}\label{lem:cgf}
$G_D$ is positive and jointly continuous: $D\times D\mapsto (0,\infty]$. 
\end{lemma}
\proof
By (\ref{eq:green:harmonic}), (\ref{eq:omega:density}), 
Lemma~\ref{lem:Hi} and symmetry, $G_D$ is locally bounded on $\{(x,y)\in D\times
D:x\neq y\}$. By Remark~\ref{r:c}, $G_D$ is locally {\it uniformly}\/
continuous in each variable, and so it is jointly continuous on this
set. Near the diagonal $\{(x,x):\,x\in D\}$ we use
\eqref{eq:green:twosets} with $U=B(x,s)\subset D$. For this $U$ the
first term on the right hand side of
\eqref{eq:green:twosets} is explicitly given by \eqref{wfg} and also
positive on $U\times U$ and the second term can be dealt with as
before. Thus, by Lemma~\ref{lem:Hi2}, $G_D(x,y)$ is jointly continuous $D\times
D\mapsto [0,\infty]$ and $G_D(x,y)>0$ on $D\times D$, regardless of
connectedness of $D$. 
\qed

For clarity we note that $G_D$ is finite and locally uniformly
continuous on $D\times D\setminus \{(x,x):\;x\in D\}$ 
(on $D\times D$ if $\alpha>d=1$), see (\ref{eq:green:twosets}) and \eqref{wfg}.

{\it Scaling}\/ will be important in what follows.
Let $k>0$. We have
$$
\int_{kU} \nu(0,y)dy=k^{-\alpha}\int_U \nu(0,y)dy\,.
$$
Similarly, if $\varphi_k(x)=\varphi(x/k)$ and $\varphi\in C^\infty_c(\Rd)$ then
$$
\Delta^{\alpha/2}\varphi_k(x)=k^{-\alpha}
\Delta^{\alpha/2}\varphi(x/k)\,,\quad x\in \Rd\,.
$$
By (\ref{eq:green:definition}) and uniqueness of the Green function we
see that
\begin{equation}
  \label{eq:sfG}
  G_{kU}(kx,kv)=k^{\alpha-d}G_U(x,v)\,,\quad x,v\in \Rd\,,
\end{equation}
hence
\begin{equation}
  \label{eq:ss}
  s_{kU}(kx)=k^\alpha s_U(x)\,,\quad x\in \Rd\,,
\end{equation}
and
\begin{equation}
  \label{eq:sjP}
  P_{kD}(kx,ky)=k^{-d}P_D(x,y)\,,\quad x,y\in \Rd\,.
\end{equation}
By (\ref{eq:green:identity}) we also have that
\begin{equation}
  \label{eq:smh}
  \omega^{kx}_{kD}(kA)=\omega^{x}_{D}(A)\,,\quad x\in \Rd,\,A\subset \Rd\,.
\end{equation}

{\it Translation invariance}\/ is equally important but easier to
observe, for example we have $G_{U+y}(x+y,v+y)=G_U(x,v)$.
Both properties enable us to reduce many of the considerations below to the setting of
the unit ball centered at the origin.


\section{Approximate factorization of Poisson kernel} \label{sec:bhp}
We keep assuming that 
{$\emptyset \neq D\subset \Rd$ is open}.  
Note that the constants in the estimates
below are independent of $D$. When $0<r\leq 1$ we denote $D_r = D \cap B_r$ and $D_r' = B^c \cup D
\setminus B_r$. 
Our first estimate is an extension of an observation made in \cite[the proof of Lemma~3.3]{bib:sw}.
\begin{lemma} \label{th:escape} 
For every $p \in (0, 1)$ there is a constant $C_{d, \alpha, p}$ such
that if $D \subset B$ then
$$
  \omega^x_D(B^c)
  \leq C_{d, \alpha, p} \, s_{D}(x)
  \, , \quad x \in D_p \,.
$$
\end{lemma}
\proof
Let $0 < p < 1$. We choose a function $\varphi \in C^\infty_c(\Rd)$
such that $0 \leq \varphi \leq 1$, 
$\varphi(y) = 1$ if $|y| \leq p$, and $\varphi(y) = 0$ if $|y| \geq 1$. Let $x \in D_p$. By \eqref{eq:green:identity} we have
\begin{eqnarray*}
  \omega^x_D(B^c)
  & = &
  \int_{B^c} (\varphi(x) - \varphi(y)) \omega^x_D(dy) \leq
  \int_{D^c} (\varphi(x) - \varphi(y)) \omega^x_D(dy)
  \\ & = &
  - \int_D G_D(x, y) \Delta^{\alpha/2} \varphi(y) dy \,.
\end{eqnarray*}
It remains to observe that $\Delta^{\alpha/2} \varphi$ is bounded and the lemma follows.
\qed

For $x \in \Rd$, $r > 0$, and a nonnegative measure $\lambda$ on $\Rd$, we let
$$
  \Lambda_{x}({\lambda}) =
  \int_{\Rd} \nu(x, y) {\lambda}(dy)
  \, , \quad \mbox{and } \quad
  \Lambda_{x, r}({\lambda}) =
  \int_{B(x, r)^c} \nu(x, y) {\lambda}(dy)\,.
$$
Note that if $k>0$ and ${\lambda}_k$ is the {\it dilation}\/ of the
measure ${\lambda}$ defined by 
\begin{equation}
  \label{eq:sk}
\int \varphi(y) {\lambda}_k(dy)={k^d}\int \varphi(ky) {\lambda}(dy)
\end{equation}
then
\begin{equation}
  \label{eq:sL}
  \Lambda_{0,kr}({\lambda}_k)=k^{-\alpha}\Lambda_{0,r}({\lambda})\,.
\end{equation}

\begin{lemma} \label{th:bounded:kernel}
Let $0<p<1$. 
There is $C_{d, \alpha, p}$ such that if 
$D \subset B$, 
$\lambda\geq 0${, $\supp \lambda\subset B^c$,} then
\begin{equation} \label{eq:bounded:kernel}
  {P}_D[\lambda](x) \leq
  C_{d, \alpha, p} \, \Lambda_{0, p}(\pim_D[\lambda])
  \, , \quad x \in D_p \,.
\end{equation}
\end{lemma}
\proof
Let $0<p<q<r\leq 1$ and $x\in D_p$.
By \eqref{eq:poisson:harmonic} and (\ref{eq:15.5}) we have
$$
  P_D[\lambda](x) =
  \Omega_{D_r}^D[P_D[\lambda], \lambda](x) \leq
  \int_{D_r'} P_{B_r}(x, y) \pim_D[\lambda](dy)
  \,.
$$
Fubini-Tonelli theorem yields
$$
  {P}_D[\lambda](x) \leq
  \frac{1}{1 - q} \int_q^1 \int_{D_r'} P_{B_r}(x,y) \pim_D[\lambda](dy) dr =
  \int_{D_q'} K(x, y) \pim_D[\lambda](dy) \, ,
$$
where, according to \eqref{eq:poisson:ball},
$$ 
  K(x, y) = \frac{1}{1 - q} \int_q^{1 \wedge |y|} P_{B_r}(x, y) dr
  = \frac{{\cal C}_{d, \alpha}}{1 - q} \int_q^{1 \wedge |y|}
    \expr{\frac{r^2 - |x|^2}{|y|^2 - r^2}}^{\alpha / 2} \frac{1}{|x - y|^d} \, dr \,.
$$
Here and below $|y|\geq q$ and $r\leq 1\wedge |y|$, which implies that
$$
  \frac{|x - y|}{|y|} \geq
  \frac{q - p}{q}
  \, , \quad
  \frac{|y| + r}{|y|} \geq 1
  \, ,\; \mbox{and } \quad
  r^2 - |x|^2 \leq 1 \,.
$$
Thus
$$
  K(x, y) \leq
  \frac{C_{d, \alpha, q / p}}{|y|^{d + \alpha / 2}} \int_q^{1 \wedge |y|}
    \frac{dr}{(|y| - r)^{\alpha / 2}} \leq
  \frac{C_{d, \alpha, q / p}}{|y|^{d + \alpha}} \,.
$$
We conclude the proof by choosing, e.g., $q = (1 + p) / 2$.
\qed

The above {\it regularization}\/ of $P_{B_r}(x, y)$ 
(\cite{bib:b:bhp}) 
is an analogue of volume averaging in classical potential theory.

\begin{lemma} \label{th:factorization}
Let $0<p < 1$. There is $C_{d, \alpha, p}$ such that 
if  $\lambda\geq 0$, {$\supp \lambda\subset B^c$} and $D \subset B$, then
\begin{equation} \label{eq:factorization}
  C_{d, \alpha, p}^{-1} \Lambda_{0, p}({\pim_D[\lambda]}) s_D(x) \leq
  {P_D[\lambda]}(x) \leq
  C_{d, \alpha, p} \Lambda_{0, p}({\pim_D[\lambda]}) s_D(x)
  \, , \quad x \in D_p\,.
\end{equation}
\end{lemma}
\proof
Let $0<p < q < r < 1$ and $x \in D_p$. By \eqref{eq:poisson:harmonic} {and \eqref{eq:n}}
we have that
\begin{equation} \label{eq:factorization:split}
  P_D[\lambda](x) = 
  \int_{D_r'} P_{D_q}(x, y) \pim_D[\lambda](dy) +
    \int_{D_r \setminus D_q} P_D[\lambda](y) \omega^x_{D_q}(dy) \,.
\end{equation}
If $v \in D_q$ and $y \in B_r^c$, then $(r - q) / r \leq |y - v| / |y|
\leq (r + q) / q$. 
Hence 
\begin{eqnarray}
{
  \int_{D_r'} P_{D_q}(x, y) \pim_D[\lambda](dy)
}
  & = & \nonumber
  \int_{D_r'} \int_{D_q}
    G_{D_q}(x, v) \nu(v, y) dv {\pim_D[\lambda](dy)}
  \\ & \asymp & \label{eq:factorization:first}
  C_{d, \alpha, r, q} \, s_{D_q}(x)
    \int_{D_r'} \nu(0, y) {\pim_D[\lambda](dy)} \,.
\end{eqnarray}
The second integral of \eqref{eq:factorization:split} is estimated
by using Lemma~\ref{th:escape},~\ref{th:bounded:kernel}, and
scaling (\ref{eq:smh},~\ref{eq:ss}):
\begin{eqnarray}
{
  \int_{D_r \setminus D_q} P_D[\lambda](y) \omega^x_{D_q}(dy)
}
  & \leq & \nonumber
  \omega^x_{D_q}(B_q^c) \sup_{D_r \setminus D_q} {P_D[\lambda]} 
  \\ & \leq & \label{eq:factorization:second}
  C_{d, \alpha, p, q, r} \, s_{D_q}(x)
    \int_{D_r'} \nu(0, y) {\pim_D[\lambda](dy)} \,.
\end{eqnarray}
Since {$P_D[\lambda]$} is nonnegative, \eqref{eq:factorization:split}, \eqref{eq:factorization:first} and \eqref{eq:factorization:second} yield:
$$
  {P_D[\lambda]}(x) \asymp
  C_{d, \alpha, p, q, r} \, s_{D_q}(x)
    \int_{D_r'} \nu(0, y) {\pim_D[\lambda](dy)} \,.
$$
Clearly, $s_{D_q}(x) \leq s_D(x)$.
In view of \eqref{eq:superharmonic}, Lemma~\ref{th:escape} and scaling we also have that
\begin{eqnarray*}
  s_D(x)
  & = &
  s_{D_q}(x) + \int_{D \setminus D_q} s_D(z) \omega^x_{D_q}(dz) \leq
  s_{D_q}(x) + \omega^x_{D_q}(B_q^c) \sup_D s_D 
  \\ & \leq &
  s_{D_q}(x) (1 + C_{d, \alpha, p, q} \sup_B s_B) =
  C_{d, \alpha, p, q} \, s_{D_q}(x)\,.
\end{eqnarray*}
Of course, $\int_{D_r'} \nu(0, y) {\pim_D[\lambda](dy)} \leq \int_{B_p^c} \nu(0, y) {\pim_D[\lambda](dy)}$. Lemma~\ref{th:bounded:kernel} yields that also
\begin{eqnarray*}
  \int_{B_p^c} \nu(0, y) {\pim_D[\lambda](dy)} & \leq &
  \int_{D_r'} \nu(0, y) {\pim_D[\lambda](dy)} +
    \frac{C_{d, \alpha} |D_r|}{p^{d + \alpha}} \sup_{D_r} {P_D[\lambda]}
  \\ & \leq &
  C_{d, \alpha, p, r} \int_{D_r'} \nu(0, y) {\pim_D[\lambda](dy)} \,.
\end{eqnarray*}
This proves
\eqref{eq:factorization} by choosing, e.g., $q=p+(1-p)/3$ and $r=p+2(1-p)/3$. 
In fact for 
{\it every} $x \in D$ we have
\begin{equation}
  \label{eq:44.5}
  {P_D[\lambda]}(x) =
  \int_{B^c} \int_D
    G_D(x, z) \nu(z, y) dz {\pim_D[\lambda](dy)} \geq
  C_{d, \alpha} \, s_D(x)
    \int_{B^c} \nu(0, y) {\pim_D[\lambda](dy)} \,.
\end{equation}
\qed
\begin{remark}\label{rem:sf}
{\rm 
{\it Scaling}\/ leaves (\ref{eq:factorization}) invariant.
Indeed, let {$\lambda_k$} be defined by
(\ref{eq:sk}) {for some $k > 0$}. By (\ref{eq:smh},~\ref{eq:sjP}),
{$P_{kD}[\lambda_k](kx) = P_D[\lambda](x)$ for $x \in D$}.
By (\ref{eq:sL}) and (\ref{eq:ss}) we have that 
$\Lambda_{0,kp}({\pim_{kD}[\lambda_k]})s_{kD}(kx)=\Lambda_{0,p}({\pim_D[\lambda]})s_D(x)$,
which is our claim. Similar observation is valid for {\it translation}.
}
\end{remark}

\begin{remark}
{\rm
The constant $C_{d,\alpha,p}$ in (\ref{eq:factorization})
may be considered nondecreasing in $p$. Indeed, if $0<p_1<p_2<1$ and
${P_D[\lambda]} \leq C_{d, \alpha, p_2} \Lambda_{0, p_2}({\pim_D[\lambda]}) s_D$ on $D_{p_2}$ then 
${P_D[\lambda]} \leq C_{d, \alpha, p_2} \Lambda_{0,
  p_1}({\pim_D[\lambda]}) s_D$ on $D_{p_1}$. 
Similarly, if
$C^{-1}_{d,\alpha,p_1}\Lambda_{0,p_1}(P^*_D[\lambda])s_D\leq
P_D[\lambda]$ on $D_{p_1}$, then 
$C^{-1}_{d,\alpha,p_2}\Lambda_{0,p_1}(P^*_D[\lambda])s_D\leq
P_D[\lambda]$ on $D_{p_1}$. 

The lower bound in (\ref{eq:factorization}) even holds with a constant
independent of $p$, see (\ref{eq:44.5}).
}  
\end{remark}


\proofof{Theorem~\ref{th:bhp:ball}} 
By (\ref{eq:sjP}) we only need to consider $r=1$.  Let $D_1 = D \cap B$. By
Lemma~\ref{th:poisson:harmonic} and (\ref{eq:omega:density}), for $i, j = 1, 2$ we have that
$P_D(x_i, y_j) = {P}_{D_1}[\lambda_j](x_i)$, where $\lambda_j\geq 0$.
Lemma~\ref{th:factorization} with $p=1/2$ yields
\begin{eqnarray*}
P_D(x_1,y_1)P_D(x_2,y_2)&\leq& 
C^2_{d,\alpha,1/2} 
\Lambda_{0,1/2}(\pim_{D_1}[\lambda_1])s_{D_1}(x_1)
\Lambda_{0,1/2}(\pim_{D_1}[\lambda_2])s_{D_1}(x_2) \\
&\leq& C^4_{d,\alpha,1/2} 
P_D(x_1,y_2)P_D(x_2,y_1)\,.\quad \qed
\end{eqnarray*}
Let $\lambda, \rho$ be nonnegative measures on $B_r^c$. Integrating
(\ref{eq:bhp}) with respect to $\lambda(dy_1)\rho(dy_2)$,
for $x_1,x_2\in D\cap B_{r/2}$ we obtain
\begin{equation}
  \label{eq:lubhp}
  P_D[\lambda](x_1)P_D[\rho](x_2)\leq C_{d,\alpha} P_D[\rho](x_1)P_D[\lambda](x_2)\,. 
\end{equation}
By translation, (\ref{eq:lubhp}) and (\ref{eq:bhp}) extend
to intersections of $D$ and balls of arbitrary center.

Inequality (\ref{eq:lubhp}) and the following {\it global}\/
version of it state our uniform BHP in a more traditional form,
see also Remark~\ref{rem:gh}.
We emphasize that the constant in (\ref{eq:tbhp})
  below does not depend on $D$, and that $D$ may be disconnected. 
\begin{theorem} \label{th:bhp}
Let $G \subset \Rd$ be open and let $K \subset G$ be compact. There is
a constant $C = C_{d, \alpha, G, K}$ with the following property. 
If $D \subset \Rd$ is open, $\lambda, \rho$ are
  nonnegative measures not charging $G\cap D^c$, 
and $f=P_D[\lambda]$, $g=P_D[\rho]$, then
\begin{equation}
  \label{eq:tbhp2}
  C^{-1} f(y)g(x)\leq
  f(x)g(y) \leq C f(y)g(x)
  \, , \quad x, y \in K\cap D\,.
\end{equation}
\end{theorem}
\proof
In what follows will use Lemma~\ref{th:factorization} and a refinement
of the argument used in the proof of the {\it global}\/
Harnack inequality (Lemma~\ref{lem:Hi2}).
For every $x \in K$ we consider a ball $B(x, r_x) \subset G$. 
We select a finite covering,
$B(x_1, p \, r_{x_1}), \dots, B(x_n, p \,
  r_{x_n})$, of $K$, where, e.g., $p = 1/2$. 
 We denote
$r_j = r_{x_j}$, $B_j = B(x_j, r_j)$, $\wt{B}_j = B(x_j, p \, r_j)$,
where $j=1,\ldots,n$, and we let $R = \diam K$ and $r = \min \{r_1,
\dots, r_n\}$. 
We now fix $x,y\in D\cap K$ and let $i,j$ be such that $x \in D \cap
\wt{B}_i$, $y \in D \cap \wt{B}_j$. Let 
$f=P_D[\lambda]$, as described above.
Note that $f$ is a Poisson integral on each $D\cap B_i$ by
(\ref{eq:poisson:harmonic}) and (\ref{eq:omega:density}).
By Lemma~\ref{th:factorization} and Remark~\ref{rem:sf}
we obtain
$$
  \frac{f(x)}{s_{D \cap B_i}(x)}
  \leq
  C_{d, \alpha, p} \expr{
    \int_{\wt{B}_i^c \setminus \wt{B}_j} \nu(x_i, z) {\pim_D[\lambda](dz)} +
    \int_{\wt{B}_i^c \cap \wt{B}_j\cap D} \nu(x_i, z) f(z) dz} \,.
$$
For $z \in \wt{B}_i^c\setminus \wt{B}_j$ we have
$|z - x_j| \leq R + |z - x_i| \leq \frac{R + rp}{rp} |z -
x_i|$, thus $\nu(x_i, z)\leq C_{d,\alpha,p,r,R}\,\nu(x_j,z)$ in the
first integral. In the second one we simply estimate
  $\nu(x_i, z) \le C_{d, \alpha,p,r}$. It follows that 
$$
  \frac{f(x)}{s_{D \cap B_i}(x)}
  \leq
  C_{d, \alpha, p, r, R} \expr{
    \int_{\wt{B}_j^c} \nu(x_j, z) {\pim_D[\lambda](dz)} +
    \int_{D \cap \wt{B}_j} f(z) d{z}}\,.
$$
We use Lemma~\ref{th:factorization} and Remark~\ref{rem:sf}
to estimate the integrals. We obtain
$$
  \frac{f(x)}{s_{D \cap B_i}(x)}
  \leq
  C_{d, \alpha, p, r, R} \expr{
    \frac{f(y)}{s_{D \cap B_j}(y)} +
    \Lambda_{x_j, p \, r_j}({\pim_D[\lambda]}) \int_{D \cap \widetilde{B}_j} s_{D \cap B_j}(z) dz}\,.
$$
But $s_{D \cap B_j}\leq s_{B_j}
\leq C_{d, \alpha,R}$, and 
$\Lambda_{x_j, p \, r_j}({\pim_D[\lambda]}) \leq C_{d, \alpha, p} \, f(y) / s_{D \cap
B_j}(y)$.
Therefore 
$$
  \frac{f(x)}{s_{D \cap B_i}(x)} \leq
  C_{d, \alpha, p, r, R} \frac{f(y)}{s_{D \cap B_j}(y)} \,.
$$
By analogous inequality for $g$ we obtain (\ref{eq:tbhp2}).
\qed
\begin{remark}\label{rem:iloraz}
{\rm
We note that (\ref{eq:tbhp2}) may be written as
\begin{equation}
  \label{eq:tbhp}
  C^{-1} \frac{f(y)}{g(y)} \leq
  \frac{f(x)}{g(x)} \leq
  C \frac{f(y)}{g(y)}
  \, , \quad x, y \in D \cap K \, ,
\end{equation}
provided the Poison integrals are nonzero and finite for (one and therefore for all)
$x\in D$. Specifically, for $f$ the condition means that
(\ref{eq:lambda}) holds and $\lambda$ is not equal to zero on $D^c$.
}
\end{remark}
\begin{remark}\label{rem:con}
{\rm 
As seen in the above proof, $C$ in  \eqref{eq:tbhp} depends only on
$d$, $\alpha$, $\diam G$ and $\dist(K, G^c)$. In fact, by scaling, its
dependence on $\diam G$ and $\dist(K, G^c)$ is only through the ratio {$\dist(K, G^c) / \diam G$}. 
}
\end{remark}
\begin{remark}\label{rem:boundedness}
{\rm Let $D \subset \Rd$ be open, $U \subset D$ bounded, 
$f ={P}_D[\lambda]$ for a nonnegative measure $\lambda$ on $D^c$. 
If $h=\dist({\rm supp}\lambda,{U})> 0$ 
and $f$ is finite at one point of $U$ then $f$ is
  {\it bounded}\/ on $U$. This follows from Theorem~\ref{th:bhp} applied to
  $K = \overline{U}$, $G=\{x\in \Rd\,:\;\dist(x,U)<h/3\}$
 and $g=P_{D \cap G}[\ind_{A}dx]\leq 1$, where
  $A=\{x\in \Rd\,:\;\dist (x,U)>2h/3\}$ is open and nonempty.
}
\end{remark}

\begin{remark}\label{rem:Greenfunction}
{\rm Let $D \subset \Rd$ be an open Greenian set, let $x_0 \in D$ and $f(x) = G_D(x,x_0)$. It is well known that the set $\{x
  \in \partial{D}: f(x) > 0\}$ is polar 
{(\cite[p.263]{bib:La})} 
so it is of Lebesgue measure zero. Let $G$ be an open bounded
  Lipschitz domain, $D \cap G \ne \emptyset$, and assume that $x_0
  \notin G$. Then $\omega_G^x(\partial{G}) = 0$ for $x \in G$. By the
  above and (\ref{eq:green:harmonic}) we have
$$
f(x) = \int_{D \setminus \overline{G}} f(w) P_{D \cap G}(x,w) \, dw, \quad x \in D \cap G,
$$
so $f$ is a Poisson integral on $D \cap G$. 
We thus may apply Theorem~\ref{th:bhp} to $D$, $G$ and $f$ as above.
We may also use Remark~\ref{rem:boundedness} for $f$. 
It follows that for arbitrary $r > 0$ the function is bounded on any
bounded subset of $D\setminus B(x_0,r)$.
The result is nontrivial if $d\leq \alpha$, especially for $d=\alpha=1$.
}  
\end{remark}


\section{Existence of limits} \label{sec:limits}

For a positive function $q$ on a nonempty set $U$ we define its relative oscillation:
$$
  \ro_U q = \ro_{x\in U}q(x)=
  \frac{\sup_{x\in U} q(x)}{\inf_{x\in U} q(x)} \,.
$$
For notational convenience, we put $\ro_U q = 1$ if $U = \emptyset$. 

The  main result of this section addresses the asymptotics of
Poisson integrals at $x=0$. (\ref{eq:nH}) gives a motivation for
(\ref{eq:oscillation}), but here $x=0$ may be, e.g., a boundary point of
$D$.
%
\begin{lemma} \label{th:oscillation}
For every $\eta > 0$ there exists $r > 0$ such that
\begin{equation} \label{eq:oscillation}
  \ro_{D \cap B_r} \frac{{P}_D[\lambda_1]}{{P}_D[\lambda_2]} \leq
  1 + \eta
\end{equation}
for all open $D \subset B$ and nonzero nonnegative measures
$\lambda_1, \lambda_2$ on $B^c$ satisfying \eqref{eq:lambda}.
\end{lemma}
\proof
Let $c$ denote $C_{d, \alpha, 1/2}$ of
Lemma~\ref{th:factorization} with $p=1/2$. 
Recall from the proof of Theorem~\ref{th:bhp:ball} 
that (\ref{eq:lubhp}) holds with $C_{d, \alpha} = c^4$. 
Thus, (\ref{eq:oscillation}) holds for $r=1/2$ with 
{$1+\eta$} 
replaced
by $c^4$. We will show that the left hand side of (\ref{eq:oscillation}) is self-improving
when $r\rightarrow 0^+$. This will be done under each of the two
complementary assumptions: (\ref{eq:inaccessible:estimate}) and
(\ref{eq:accessible:estimate}) below. First, however, we need some preparation.
For $0 < p < q < 1/2$ and a measure ${\lambda}$ let $D_{p,q}=D_q \setminus
D_p$ and
$$
\Lambda_{x, p, q}({\lambda})= \int_{D_{p, q}} \nu(x, y) {\lambda}(dy)\,.
$$
{We also denote}
$$
  \begin{array}{ccc}
  f_i = P_D[\lambda_i] \, , \quad &
  f_i^{p, q} = P_{D_p}[\ind_{D_{p, q}}\pim_D[\lambda_i]] \, , \quad &
  \wt{f}_i^{p, q} = P_{D_p}[\ind_{D_q'}\pim_D[\lambda_i]] \, , \\
  f_i^* = \pim_D[\lambda_i] \, , \quad &
  f_i^{p, q*} = \pim_{D_p}[\ind_{D_{p, q}}\pim_D[\lambda_i]] \, , \quad &
  \wt{f}_i^{p, q*} = \pim_{D_p}[\ind_{D_q'}\pim_D[\lambda_i]] \,.
  \end{array}
$$
What follows will be valid for $i=1$ and $i=2$. 
By \eqref{eq:poisson:harmonic} and (\ref{eq:omega:density}) we have
$f_i = f_i^{p, q} + \wt{f}_i^{p, q}$
and $f^*_i = f_i^{p, q*} + \wt{f}_i^{p, q*}$.
For $r \in (0, 1/2]$ we denote $m_r = \inf_{D_r} (f_1 / f_2)$ and $M_r
= \sup_{D_r} (f_1 / f_2)$. As we noted above, 
$M_r \leq c^4 m_r$.
Let $\varepsilon > 0$.

Let $q \in (0, 1/2]$ and let $p = p(q) \in (0, q /2)$ (depending on $q$
and $\varepsilon$) be defined by 
\begin{equation}
  \label{eq:wq}
(q + 2 p) / (q - 2 p) = 1 + \varepsilon\,,
\end{equation}
so that if  $z \in D_{2 p}$ and $y \in B^c_q$ then $(1 + \varepsilon)^{-d - \alpha}
\nu(0, y) \leq \nu(z, y) \leq (1 + \varepsilon)^{d + \alpha} \nu(0,
y)$. Thus, for $x \in D_{2 p}$ we have
$$
  \wt{f}_i^{2 p, q}(x) =
  \int_{D_q'} \int_{D_{2 p}} G_{D_{2 p}}(x, z) \nu(z, y) dz {f_i^*(dy)} \leq
  (1 + \varepsilon)^{d + \alpha} \Lambda_{0, q}(f_i{^*}) s_{D_{2 p}}(x) \, ,
$$
and 
$$
  \wt{f}_i^{2 p, q}(x) \geq
  (1 + \varepsilon)^{-d - \alpha} \Lambda_{0, q}(f_i{^*}) s_{D_{2 p}}(x) \,.
$$

We will now examine consequences of the following assumption:
\begin{equation} \label{eq:inaccessible:estimate}
  \Lambda_{0, p, q}(f_i^*) \leq
  \varepsilon \, \Lambda_{0, q}(f_i^*) \,,\quad i=1,2\,.
\end{equation}
If \eqref{eq:inaccessible:estimate} holds then using 
Lemma~\ref{th:factorization}
and Remark~\ref{rem:sf} we obtain
$$
  f_i^{2 p, q}(x) \leq
  c \, s_{D_{2 p}}(x) \Lambda_{0, p}(f_i^{2 p, q{*}}) \leq
  c \, s_{D_{2 p}}(x) \Lambda_{0, p, q}(f_i{^*}) \leq
  c \, \varepsilon \, s_{D_{2 p}}(x) \Lambda_{0, q}(f_i^*)
  \, , \quad x \in D_p \,.
$$
Recall that $f_i = f_i^{2 p, q} + \wt{f}_i^{2 p, q}$. 
Thus, if \eqref{eq:inaccessible:estimate} holds then we have
\begin{equation} \label{eq:inaccessible:explicit:estimate}
  \frac{(1 + \varepsilon)^{-d - \alpha} \Lambda_{0, q}(f_1{^*})}
    {(c \, \varepsilon + (1 + \varepsilon)^{d + \alpha}) \Lambda_{0, q}(f_2{^*})} \leq
  \frac{f_1(x)}{f_2(x)} \leq
  \frac{(c \, \varepsilon + (1 + \varepsilon)^{d + \alpha}) \Lambda_{0, q}(f_1{^*})}
    {(1 + \varepsilon)^{-d - \alpha} \Lambda_{0, q}(f_2{^*})}
  \, , \quad x \in D_p \, ,
\end{equation}
and, finally,
\begin{equation} \label{eq:inaccessible:oscillation}
  \ro_{D_p} \frac{f_1}{f_2} \leq 
  (c \, \varepsilon + (1 + \varepsilon)^{d + \alpha})^2 (1 + \varepsilon)^{2 d + 2 \alpha}
  \,.
\end{equation}
We are satisfied with (\ref{eq:inaccessible:oscillation}) for the moment.
Let $0 < \bar{p} < \bar{q}/4<\bar{q} < 1/2$, $g = f_1^{2 \bar{p},
  \bar{q}} - m_{\bar{q}} f_2^{2 \bar{p}, \bar{q}}$, and $h = M_{\bar{q}} f_2^{2 \bar{p}, \bar{q}} - f_1^{2
  \bar{p}, \bar{q}}$. Note that on $D_{2 \bar{p}}$ both $g$ and $h$ are Poisson integrals 
of nonnegative measures. 
If $D_{\bar{p}} \neq \emptyset$ then by Theorem~\ref{th:bhp:ball} (see
also (\ref{eq:lubhp})),
$$
  \sup_{D_{\bar{p}}} \frac{f_1^{2 \bar{p}, \bar{q}}}{f_2^{2 \bar{p}, \bar{q}}} - m_{\bar{q}} =
  \sup_{D_{\bar{p}}} \frac{g}{f_2^{2 \bar{p}, \bar{q}}} \leq
  c^4 \inf_{D_{\bar{p}}} \frac{g}{f_2^{2 \bar{p}, \bar{q}}} = 
  c^4 \expr{\inf_{D_{\bar{p}}} \frac{f_1^{2 \bar{p}, \bar{q}}}{f_2^{2 \bar{p}, \bar{q}}} - m_{\bar{q}}} \, ,
$$
and
$$
  M_{\bar{q}} - \inf_{D_{\bar{p}}} \frac{f_1^{2 \bar{p}, \bar{q}}}{f_2^{2 \bar{p}, \bar{q}}} =
  \sup_{D_{\bar{p}}} \frac{h}{f_2^{2 \bar{p}, \bar{q}}} \leq
  c^4 \inf_{D_{\bar{p}}} \frac{h}{f_2^{2 \bar{p}, \bar{q}}} = 
  c^4 \expr{M_{\bar{q}} - \sup_{D_{\bar{p}}} \frac{f_1^{2 \bar{p}, \bar{q}}}{f_2^{2 \bar{p}, \bar{q}}}} \,.
$$
By adding these inequalities we obtain
\begin{equation} \label{eq:accessible:base}
  (c^4 + 1)\expr{\sup_{D_{\bar{p}}} \frac{f_1^{2 \bar{p}, \bar{q}}}{f_2^{2 \bar{p}, \bar{q}}} -
    \inf_{D_{\bar{p}}} \frac{f_1^{2 \bar{p}, \bar{q}}}{f_2^{2 \bar{p}, \bar{q}}}} \leq
  (c^4 - 1)(M_{\bar{q}} - m_{\bar{q}}) \,.
\end{equation}

We will now examine consequences of the following assumption:
\begin{equation} \label{eq:accessible:estimate}
  \Lambda_{0, \bar{q}}(f_i^*) \leq
  \varepsilon \Lambda_{0, 2 \bar{p}, \bar{q}}(f_i^*) \,.
\end{equation}
By Lemma~\ref{th:factorization} and Remark~\ref{rem:sf}
(consider $D=D_{2\bar{p}}$, $p=2\bar{p}$ and $B_{\bar{q}}$ replacing
$B$ therein),
for $x \in D_{\bar{p}}$ we have
$$
  \wt{f}_i^{2 \bar{p}, \bar{q}}(x) 
\leq c \, s_{D_{2 \bar{p}}}(x) \Lambda_{0, 2\bar{p}}(\wt{f}_i^{2 \bar{p}, \bar{q}{*}}) 
= c \, s_{D_{2 \bar{p}}}(x) \Lambda_{0, \bar{q} / 2}(\wt{f}_i^{2 \bar{p}, \bar{q}{*}}) 
\leq c \, s_{D_{2 \bar{p}}}(x) \Lambda_{0, \bar{q}}(f_i{^*})\,.
$$
Hence, by our assumption (\ref{eq:accessible:estimate}),
Lemma~\ref{th:factorization} and Remark~\ref{rem:sf} applied to
$D_{\bar{p}}\subset B_{2\bar{p}}$
$$
  \wt{f}_i^{2 \bar{p}, \bar{q}}(x) \leq
  c \, \varepsilon \, s_{D_{2 \bar{p}}}(x) \Lambda_{0, 2 \bar{p}, \bar{q}}(f_i{^*}) \leq
  c \, \varepsilon \, s_{D_{2 \bar{p}}}(x) \Lambda_{0, \bar{p}}(f_i^{2 \bar{p}, \bar{q}{*}}) \leq
  c^2 \varepsilon \, f_i^{2 \bar{p}, \bar{q}}(x)
  \, , \quad x \in D_{\bar{p}} \,.
$$
Since $f_i = f_i^{2 \bar{p},\bar{q}} + \wt{f}_i^{2 \bar{p}, \bar{q}}$
on $D_{2\bar{p}}$, this and
\eqref{eq:accessible:base} yield
$$
(c^4 + 1) \left(M_{\bar{p}}/(1 + c^2 \varepsilon)-
m_{\bar{p}}(1 + c^2 \varepsilon)\right)
\leq
(c^4 - 1) (M_{\bar{q}} - m_{\bar{q}}) \,.
$$
Note that $m_{\bar{p}} \geq m_{\bar{q}}$. Dividing by $m_{\bar{q}}$ finally gives
\begin{equation} \label{eq:accessible:oscillation}
  \ro_{D_{\bar{p}}} \frac{f_1}{f_2} \leq
  (1 + c^2 \varepsilon)^2 +
    (1 + c^2 \varepsilon) 
\frac{c^4 - 1}{c^4 + 1}
    \expr{\ro_{D_{\bar{q}}} \frac{f_1}{f_2} - 1} \,.
\end{equation}

We now come to the conclusion of our considerations.
Let $\eta > 0$. 
If $\varepsilon$ is small enough then the right hand side of \eqref{eq:inaccessible:oscillation} is smaller
than $1 + \eta$ and right hand side of
\eqref{eq:accessible:oscillation} does not exceed $\varphi(\ro_{D_{\bar{q}}}
(f_1 / f_2))$, where
$$
  \varphi(t) = 1 + \frac{\eta}{2} + \frac{c^4}{c^4 + 1} (t - 1)
  \,,\quad t\geq 1\,.
$$
Let $\varphi^{1}=\varphi$, $\varphi^{l+1}=\varphi(\varphi^{l})$, $l=1,2,\ldots$.
Observe that $\varphi(t)=t$ for 
$t=1+\eta(c^4+1)/2$, and $\varphi(t)<t$ for $t>1+\eta(c^4+1)/2$. 
Thus the $l$-fold compositions $\varphi^l(c^4)$ converge to $1+\eta(c^4+1)/2$ as
$l\rightarrow \infty$. In what follows let $l$ be such that $$\varphi^l(c^4)<1+\eta(c^4+1)\,.$$
Let $k$ be the least integer such that $k - 1 >
c^2/\varepsilon^{2}$. We denote $n = lk$. 
Let $q_0 = 1/2$, $q_{j + 1}= p(q_j)$ for $j=0,\ldots,n-1$ (see
(\ref{eq:wq})), and $r = q_n$. If for any $j < n$,
\eqref{eq:inaccessible:estimate} holds with $q = q_j$ and $p = p(q) = q_{j +
  1}$, then
$$
  \ro_{D_r} \frac{f_1}{f_2} \leq
  \ro_{D_{q_{j+1}}} \frac{f_1}{f_2} \leq
  1 + \eta \,,
$$
and we are done by the definition of $\varepsilon$ and (\ref{eq:inaccessible:explicit:estimate}).
Otherwise for $j=0,\ldots, n-1$, we have $\Lambda_{0, q_{j +
    1}, q_j}(f_i{^*}) > \varepsilon \Lambda_{0, q_j}(f_i{^*})$ for $i = 1$ or
$i = 2$. Note that by Lemma~\ref{th:factorization}
$$
  c^{-1}\frac{f_i(x)}{\Lambda_{0, q_j}(f_i{^*})} \leq
  s_{D_{2 q_j}}(x) \leq
  c \frac{f_{3-i}(x)}{\Lambda_{0, q_j}(f_{3-i}{^*})}
  \, , \quad x \in D_{q_{j + 1}, q_j} \,.
$$
Hence 
$
\Lambda_{0, q_{j+1},q_j}(f_i{^*})/
\Lambda_{0,q_j}(f_i{^*})
\leq
c^2 \Lambda_{0, q_{j+1},q_j}(f_{3-i}{^*})/
\Lambda_{0,q_j}(f_{3-i}{^*}) 
$,
and so
$\Lambda_{0, q_{j + 1}, q_j}(f_i{^*}) > c^{-2} \varepsilon \,
\Lambda_{0, q_j}(f_i{^*})$ for both $i = 1$ and $i=2$ (and all
$j=0,1,\ldots, n-1$). 
If $0\leq j<l$ and $\bar{p} = q_{(j + 1) k}$, $\bar{q} = q_{j k}$, then
$$
  \Lambda_{0, 2 \bar{p}, \bar{q}}(f_i{^*}) \geq
  \Lambda_{0, q_{(j + 1) k - 1}, q_{j k}}(f_i{^*}) \geq
  (k - 1) c^{-2} \varepsilon \, \Lambda_{0, \bar{q}}(f_i{^*}) \geq
  \varepsilon^{-1} \Lambda_{0, \bar{q}}(f_i{^*}) \, ,
$$
so that \eqref{eq:accessible:estimate} is satisfied. We conclude that
\eqref{eq:accessible:oscillation} holds. 
Recall that $\ro_{D_{1/2}} (f_1 / f_2) \leq c^4$. By the definition of
$l$ and monotonicity of $\varphi$
$$
  \ro_{D_{q_{lk}}} \frac{f_1}{f_2} \leq
  \varphi \expr{\ro_{D_{q_{(l - 1)k}}} \frac{f_1}{f_2}} \leq
  \ldots \leq
  \varphi^l \expr{\ro_{D_{q_0}} \frac{f_1}{f_2}} \leq
  1 + \eta (c^4+1) \, ,
$$
i.e. $\ro_{D_r} (f_1 / f_2) \leq 1 + \eta(c^4+1)$. Since $\eta>0$ was
arbitrary, the proof is complete.
\qed

\begin{corollary} \label{th:limit}
If $D$ is bounded, $0 \in \partial D$ and $y \in D^c\setminus\{0\}$, then
\begin{equation}\label{eq:limit}
  \lim_{D \ni x \rightarrow 0} \frac{P_D(x, y)}{s_D(x)}\;\;\mbox{exists}.
\end{equation}
\end{corollary}
\proof
We may assume that $P_D(x,y)<\infty$ for $x\in D$.
For bounded $D$, by (\ref{eq:poisson:definition}) and (\ref{eq:ds}) we
have that  
\begin{equation}
  \label{eq:ws}
s_D(x)=\lim_{|z|\rightarrow
  \infty}P_D(x,z)/\nu(0,z)\,. 
\end{equation}
In fact, $P_D(x,z)/(s_D(x)\nu(0,z))\to 1$ uniformly in $x\in D$.
{By scaling, we may assume that $|y| \ge 1$.}
We apply Lemma~\ref{th:oscillation} to $\lambda_1 = \varepsilon_{y}$
and $\lambda_2 = \varepsilon_{z}/\nu(0,z)$. It follows that $\ro_{D_r} P_D(\cdot,y) / s_D(\cdot) \rightarrow 1$
as $r \rightarrow 0^+$, which, in presence of (\ref{eq:factorization}), is equivalent to the convergence 
to a finite, positive limit.
\qed

As an addition to Corollary~\ref{th:limit} we note that if $0$ is inaccessible
from $D$, then we have
\begin{equation} \label{eq:limit:explicite}
  \lim_{D \ni x \rightarrow 0} \frac{{P_D[\lambda_1]}(x)}{{P_D[\lambda_2]}(x)} =
  \frac{\int_{\Rd} \nu(0, y) {\pim_D[\lambda_1](dy)}}{\int_{\Rd} \nu(0, y) {\pim_D[\lambda_2](dy)}}\,,
\end{equation}
and 
\begin{equation}
  \label{eq:gins}
  \lim_{D \ni x \rightarrow 0} \frac{{P_D[\lambda_1]}(x)}{s_D(x)} =
  \frac{\int_{\Rd} \nu(0, y) {\pim_D[\lambda_1](dy)}}{\int_{\Rd} \nu(0, y) s_D(y) dy + 1}\,.
\end{equation}
Here {$\lambda_1$, $\lambda_2$ are nonnegative measures on $B^c$} for
which the Poisson integrals are positive and finite. 
Indeed, by Lemma~\ref{th:factorization} the integrals $\int \nu(0, y)
{\pim_D[\lambda_i](dy)}$ are finite. Hence, for every $\varepsilon >
0$ we can find $q > 0$ such that
(\ref{eq:inaccessible:estimate}) is satisfied with $p = q/2$. It follows that (\ref{eq:inaccessible:explicit:estimate})
holds. Since $\varepsilon$ was arbitrary, the first equality is
proved. The second one follows by using (\ref{eq:ws}).

We remark in passing that (\ref{eq:ws}) yields (\ref{eq:sBr})
as a consequence of (\ref{eq:poisson:ball}), and sheds
some light on the role of $s_D$ as a substitute, at infinity, for the Poisson kernel, see Theorem~\ref{th:martin}.

\section{Harmonicity}\label{sec:harmonicity}
\begin{definition}\label{def:harm}
Let $f$ be a nonnegative continuous function on $D$ and
let $\lambda$ be a nonnegative measure on $D^c$.
We say that $f$ is $\alpha$-harmonic in $D$ with outer
charge $\lambda$ on $D^c$ if for every open bounded $U$ such that
$\overline{U}\subset D$
we have (see (\ref{eq:cm}))
\begin{equation} \label{eq:harmonic:definition}
  f(x) ={\Omega}^D_U[f,\lambda](x) , \quad x \in U\,.
\end{equation}
\end{definition}

We note that the integral in (\ref{eq:harmonic:definition})
is finite by the {\it assumption}\/ that the left  hand side of
(\ref{eq:harmonic:definition}) is continuous (hence finite).
Thus for (nonnegative) $f$ which is $\alpha$-harmonic on
(nonempty open) $D$ with outer charge $\lambda\geq 0$, by \eqref{eq:poisson:ball} and
\eqref{eq:omega:density} we necessarily have that
\begin{equation}
  \label{eq:calkfh}
\int_D (1+|y|)^{-d-\alpha}f(y)dy
+\int_{D^c}(1+|y|)^{-d-\alpha}\lambda(dy)<\infty\,.
\end{equation}

The present definition extends the usual definition of an $\alpha$-harmonic
function (\cite{bib:b:bhp}) 
by allowing measures as ``boundary values'' (outer charge) on $D^c$. 
The (genuine) $\alpha$-harmonic functions studied so far in the
  literature  correspond to absolutely continuous measures $\lambda$.
For such measures we can {\it denote} $d \lambda/d x$ by $f$ on $D^c$, and
(\ref{eq:harmonic:definition}) then reads
\begin{equation} \label{eq:harmonic:definition-ac}
  f(x) =\int_{U^c}f(y)\omega^x(dy)\, , \quad x \in U\,,
\end{equation}
for open precompact $U\subset D$, see (\ref{eq:n}).
We consider (\ref{eq:harmonic:definition}),
(\ref{eq:harmonic:definition-ac}) a mean-value property
because $\omega^x_U$ is a probability measure.

Formula~(\ref{eq:green:harmonic}) yields that the function $x\mapsto G_D(x, y)$
is $\alpha$-harmonic on $D \setminus \{y\}$ with zero outer charge,
see also Remark~\ref{r:c}.
Also, $x\mapsto\omega^x_D(A)$ is $\alpha$-harmonic on $D$ for every
set $A$ by \eqref{eq:omega:harmonic}, and (\ref{eq:n}) applied to
$U$. Here the outer charge is $\ind_A dx$.  
Lastly, 
if $f=P_D(\cdot,y)$ is finite on $D$ then 
by (\ref{eq:harm}) it
is $\alpha$-harmonic in $D$ with outer charge $\varepsilon_y$. 
If a measure $\lambda$ is nonnegative and 
$P_D[\lambda]$ is
finite on $D$ then by Fubini-Tonelli and Remark~\ref{r:c},
$P_D[\lambda]$ is $\alpha$-harmonic in $D$ with outer
charge $\lambda$.
By the same token, 
$x \mapsto P_D(x, y)$ is {\it not} $\alpha$-harmonic in $D$
with zero outer charge. 
To be absolutely clear on this, we consider $D = B$ and $|y|>1$. We note that the function
$x\mapsto P_B(x, y)$ vanishes on $B^c$, and has a maximum inside $B$.
Thus the mean-value property \eqref{eq:harmonic:definition} 
cannot hold with $\lambda=0$ (recall that $\omega^x_U(\Rd)=1$). 
As we will see below, not every harmonic function is a
Poisson integral, neither is every harmonic function on $D$ represented by an
integral against the harmonic measure of $D$.

If $f$ is a function on $\Rd$ continuous on $D$ and
$\lambda(dx) = f(x) dx$ on $D^c$, then
\eqref{eq:harmonic:definition} is equivalent to 
\begin{equation}
  \label{eq:dh}
\Delta^{\alpha / 2} f(x) = 0\,,\quad x \in D\,.
\end{equation}
The result is given in \cite{BB1}, and its proof can be extended 
to the present more general setting. However, we will not use (\ref{eq:dh})
in the sequel, and we leave the verification of the extension
to the interested reader. 
We also refer the reader to \cite{bib:bs} to see
the limitations of {\it pointwise}\/ definition of harmonic 
functions by means of the corresponding generators.


%

\begin{remark}\label{rem:har}
{\rm
If open precompact $D'\subset D$ is a Lipschitz domain, then
$\omega^x_{D'}(\partial D')=0$ (see the discussion following
(\ref{eq:n})), and $f$ in Definition~\ref{def:harm} can be considered a Poisson 
integral on $D'$. This immediately yields Harnack inequality for general
(nonnegative) harmonic functions, see (\ref{eq:nH}), (\ref{eq:nH2}), and Remark~\ref{r:c}. 
}  
\end{remark}

In a number of considerations below we need to allow $U$ touching $\partial D$
in (\ref{eq:harmonic:definition}). 
Let $D^{(r)}$ be the set of {\it regular}\/ boundary points for $D$ (see
Section~\ref{sec:preliminaries}). 
It is known that $\omega^x_U(\partial D\setminus D^{(r)}) = 0$
for every open $U$ and $x\in U$, 
and $|\partial D \setminus D^{(r)}| = 0$
(\cite{bib:La}).

\begin{lemma} \label{th:omega:twofunctions}
Suppose that $0 \leq g \leq f$ on $D$,
and $f, g$ are $\alpha$-harmonic on $D$ with zero outer
charge. If $U \subset D$ and 
$f(x) = \int_{U^c} f(y) \omega^x_U(dy)$, $x \in U$, then 
$g(x) = \int_{U^c} g(y)\omega^x_U(dy)$, $x \in U$.
\end{lemma}
\proof
Let $D_n$ be an increasing sequence of open sets precompact in 
$D$ such that $D = \bigcup D_n$. Let $U_n = U \cap D_n$,
so that $U_n$ increase to (possibly unbounded) $U$. Then
$$
  0 \leq
  \int_{U \setminus U_n} g(y) \omega^x_{U_n}(dy) \leq
  \int_{U \setminus U_n} f(y) \omega^x_{U_n}(dy) =
  f(x) - \int_{D \setminus U} f(y) \omega^x_{U_n}(dy)
  \, , \quad x \in U_n \,.
$$
Recall that $\omega^x_{U_n}$ increase on $D\setminus U$ and vaguely
converge to $\omega^x_U$ on $\Rd$. Since $D\setminus U \subset D$ 
and $\bigcap \cl{U \setminus U_n} \subset \partial D$ are disjoint, we conclude that
actually $\omega^x_{U_n}$ increase to $\omega^x_U$ on $D \setminus U$.
%
%
Hence,
by monotone convergence, the right hand side tends to $0$ when $n \rightarrow
\infty$. Thus
$$
  g(x) =
  \lim_{n \rightarrow \infty} \int_{D \setminus U_n} g(y) \omega^x_{U_n}(dy) =
  \lim_{n \rightarrow \infty} \int_{D \setminus U} g(y) \omega^x_{U_n}(dy) =
  \int_{D \setminus U} g(y) \omega^x_U(dy) \,.
  \qed
$$
%
%
\begin{lemma} \label{th:omega:twosets}
Let $D_1$, $D_2$ be open sets such that 
\begin{equation}
  \label{eq:ss}
\dist (D_1\setminus D_2\,,\;D_2\setminus D_1)>0\,.
\end{equation}
Let $D=D_1\cup D_2$ and assume that $\omega^x_D(D^c)>0$
for (one and therefore for all) $x \in D$.
Let $f\geq 0$ be a function on $\Rd$ such that
$f=0$ on $D^c$, and for $i=1,2$ we have
\begin{equation}
  \label{eq:wsl}
f(x) =\int f(y) \omega^x_{D_i}(dy)\quad \mbox{if }\; x \in
D_i\,.
\end{equation}
If $D_1$ is bounded then $f=0$ on the whole of $D$.
\end{lemma}
\noindent
\proof
Note that effectively the integration in (\ref{eq:wsl}) is only over
$D_{3-i}\setminus D_i$, $i=1,2$. We can assume that both these sets
are nonempty.
We then observe that $f$ is {\it bounded}\/ on $D_1\setminus D_2$. 
Indeed, if we consider {\it Lipschitz}\/ open set $U$ such that $\overline{D_1\setminus D_2}\subset
U\subset D_1\cup D^c$ and $\dist(U,D_2\setminus D_1)>0$ then it follows from
(\ref{eq:wsl}) with $i=1$, (\ref{eq:omega:harmonic}) and (\ref{eq:omega:density}) that $f$ is a
Poisson integral on $D_1\cap U$, and the boundedness follows from Remark~\ref{rem:boundedness}.   
By considering $i=2$ in (\ref{eq:wsl}) we see that $f$ is bounded on
$D_2$, hence on the whole of $D$.

Let $\{X_t\,,\,t\geq 0;\,\pr_x\,,\, \ex_x\,,\,x\in \Rd\}$ be the
isotropic stable L\'evy process with the
corresponding Markov probabilities and expectations. 
For open set $U$ we define
the first {\it entrance} time of $U^c$, $\tau_U=\inf\{t\geq 0\,:\;X_t\notin U\}$. 
It is well known that
$\omega^x_{U}(A)=\pr_x(\tau_U<\infty\,,\,X_{\tau_U}\in A)$ for $x\in
U$.
Let $T_1=\tau_{D_1}$. We define ``shuttle'' times
\begin{eqnarray*}
T_{2n}&=&T_{2n-1}+\tau_{D_2}\circ \Theta_{T_{2n-1}}\,,\\
T_{2n+1}&=&T_{2n}+\tau_{D_1}\circ \Theta_{T_{2n}}\,,\quad n=1,2,\ldots\,,
\end{eqnarray*}
where $\Theta$ is the usual shift operator: $(X\circ\Theta_{s})_t=X_{s+t}$.
If $X_{T_{n}}\in D^c$ then 
\begin{equation}
  \label{eq:fn}
\tau_D=T_n=T_{n+1}=T_{n+2}=\ldots\,,\quad \mbox{and } \quad
X_{\tau_D}=X_{T_n}=X_{T_{n+1}}=\ldots.
\end{equation}
Otherwise, 
\begin{equation}
  \label{eq:ncc}
  |X_{T_{n+1}}-X_{T_n}|\geq
\dist (\overline{D_1\setminus D_2},\overline{D_2\setminus D_1})
=\dist (D_1\setminus D_2,D_2\setminus D_1)>0\,.
\end{equation}
We define Markov time $T_\infty=\lim_{n\to \infty}T_n$.
Clearly, $T_\infty\leq\tau_D$.
By quasi-left-continuity of $\{X_t\}$ we have that 
$X_{T_n}\to X_{T_\infty}$ as $n\to \infty$, if $T_\infty <\infty$.
But (\ref{eq:ncc}) then implies (\ref{eq:fn}), in particular
$T_\infty=\tau_D$. Clearly, if $T_\infty=\infty$, in particular if
$T_n=\infty$ for some $n$, then we also have $T_\infty=\tau_D$.
If $T_{n-1}=T_n<\infty$ then $X_{T_n}\in D^c$, and so
by strong Markov property of $\{X_t\}$ and (\ref{eq:wsl}) we obtain
\begin{eqnarray}
  f(x)&=&\ex_x \left\{
T_n<\infty
\,;\;
f(X_{T_n})
\right\}
=\ex_x \left\{
T_{n-1}<T_n<\tau_D
\,;\;
f(X_{T_n})
\right\}\nonumber\\
&\leq& \|f\|_\infty\, \pr_x\left\{T_{n-1}<T_n<\tau_D\right\}\,,\quad
n\geq 2\,.\label{eq:ofg}
\end{eqnarray}
Here $\|f\|_\infty=\sup_{x\in \Rd} f(x)$. We consider
\begin{eqnarray}
g_n(x)&=&\pr_x\left\{T_{n-1}<T_n<\tau_D\right\}
=\pr_x\left\{T_{n-3}<T_{n-2}<T_{n-1}<T_n<\tau_D\right\}\nonumber\\
&=&\ex_x\left\{
T_{n-3}<T_{n-2}<\tau_D
\,;\;
\pr_{X_{T_{n-2}}}\left[T_1<T_2<\tau_D\right]
\right\}\nonumber\\
&=&\ex_x\left\{
T_{n-3}<T_{n-2}<\tau_D
\,;\;
g_2(X_{T_{n-2}})
\right\}
\leq  \|g_{n-2}\|_\infty\,\|g_2\|_\infty\,, \label{eq:ig2}
\end{eqnarray}
where $n\geq 4$. We will examine
\begin{eqnarray}
g_2(x)&=&\pr_x\left\{T_{1}<T_2<\tau_D\right\} \nonumber\\
&=&\ex_x\left\{X_{\tau_{D_1}}\in
D_2\,;\;\pr_{X_{\tau_{D_1}}}[\tau_{D_2}<\infty\,,\,X_{\tau_{D_2}}\in D_1]
\right\}\,.
  \label{eq:rg2}
\end{eqnarray}
If $g_2\equiv 1$ on $D$ then $g_n\equiv 1$ on $D$ for every even
$n$, see (\ref{eq:ig2}). By the definition of $g_n$ and our discussion
of (\ref{eq:ncc}), we then have that $\pr_x\left\{\tau_D=\infty\right\}=1$, 
contradicting the assumptions of the lemma.
Thus, $g_2(x)<1$ for some $x\in D$.
Observing the expressions in 
(\ref{eq:rg2}), by Harnack inequality we conclude that
\begin{equation}
  \label{eq:c1}
  \pr_x\left\{X_{\tau_{D_1}}\in D_2\right\}<1\quad \mbox{for all $x\in D_1$,}
\end{equation}
or
\begin{equation}
  \label{eq:c2}
  \pr_x\left\{\tau_{D_2}<\infty\,,\, X_{\tau_{D_2}}\in D_1\right\}<1\quad \mbox{for all $x\in D_2$.}
\end{equation}

Assume that (\ref{eq:c1}) holds. By Theorem~\ref{th:bhp} applied to
the above mentioned set $U\cap D_1$, there is $c>0$ such that
$$
\pr_x\left\{X_{\tau_{D_1}}\in D_2\right\}
\leq c
\pr_x\left\{X_{\tau_{D_1}}\notin
  D_2\right\}\,,\quad x\in D_1\setminus D_2\,,
$$
thus 
\begin{equation}
  \label{eq:c1o}
\pr_x\left\{X_{\tau_{D_1}}\in D_2\right\}\leq 1-c\,,\quad x\in D_1\setminus D_2\,,
\end{equation}
and so $\|g_2\|_\infty<1-c$.

Assuming (\ref{eq:c2}), there will be a ball $B'\subset D_2\setminus D_1$ such that 
$\ex\left\{X_{\tau_{D_2}}\in D_1\right\}\leq 1-\varepsilon$ for some $\varepsilon>0$.
If (\ref{eq:c1}) is false then
$\pr_x\left\{X_{\tau_{D_1}}\in D_2\right\}=1=\pr_x\left\{X_{\tau_{D_1}}\in D_2\setminus D_1\right\}
=\pr_x\left\{X_{\tau_{D_1\cup D^c}}\in D_2\setminus D_1\right\}$
for all $x\in D_1$.
By (\ref{eq:poisson:ball}) there is $c>0$, independent of $x\in
D_1\setminus D_2$, such that 
$\pr_x\left\{X_{\tau_{D_1}}\in B'\right\}\geq
\pr_x\left\{X_{\tau_{B(x,r)}}\in B'\right\}\geq c$, where
$r=\dist(x,D_2\setminus D_1)$.
We obtain that $\|g_2\|_\infty\leq 1-c\varepsilon$.

Thus, (\ref{eq:c1}) and (\ref{eq:c2}) imply that $\|g_{2n}\|_\infty\to
0$, and so $f\equiv 0$, see (\ref{eq:ofg}), (\ref{eq:ig2}).
\qed

Lemma~\ref{th:omega:twosets} applies, e.g., if $D_1$, $D_2$ are
overlapping finite open intervals on the line.

\begin{lemma} \label{th:regular:harmonic}
Let $f$ be a nonnegative function on $\Rd$ which is
$\alpha$-harmonic in bounded $D$, with outer charge $f(x) dx$.
Suppose that $f$ is
bounded and continuous on $D\cup D^{(r)}$. Then 
\begin{equation}
  \label{eq:rh}
f(x) = \int f(y) \omega^x_D(dy)\,,\quad x\in D\,.  
\end{equation}
\end{lemma}
\proof
Let $D_n$ be an increasing sequence of open sets precompact in $D$
such that $\bigcup D_n = D$. Recall that
$\omega^x_{D_n}\to\omega^x_D$ weakly as $n\to \infty$ {and $\omega^x_D(\cl{D} \setminus D^{(r)}) = 0$}.
Since ${f}$ is bounded and continuous on 
{$D \cup D^{(r)}$},
$$
  \lim\limits_{n \rightarrow \infty} \int\limits_{\cl{D}} f(y) \omega^x_{D_n}(dy) =
  \lim\limits_{n \rightarrow \infty} \int\limits_{D\cup D^{(r)}} f(y) \omega^x_{D_n}(dy) =
  \int\limits_{D\cup D^{(r)}} f(y) \omega^x_D(dy)=
  \int\limits_{\cl{D}} f(y) \omega^x_D(dy)\,.
$$
%
%
We obtain (\ref{eq:rh}) by noting that $P_{D_n}(x, y) \nearrow P_D(x,
y)$ for $x \in D$, $y \in D^c$, and so
$$
  \lim\limits_{n \rightarrow \infty} \int\limits_{\cl{D}^c} P_{D_n}(x, y) f(y) dy =
  \int\limits_{\cl{D}^c} P_D(x, y) f(y) dy\,.
  \qed
$$

We remark in passing that (\ref{eq:rh}) implies $\alpha$-harmonicity through
(\ref{eq:omega:harmonic}). The reverse implication is not true as we
will see from the example of the Martin kernel with the pole at an accessible
boundary point.


\section{Martin kernel} \label{sec:martin}
It follows from \eqref{eq:calkfh}, (\ref{eq:green:harmonic}) and (\ref{eq:n}) that for open Greenian $D$,  
\begin{equation}
  \label{eq:gc}
\int_\Rd G_D(x,v)(1+|v|)^{-d-\alpha}dv<\infty\,,\quad x\in \Rd\,.
\end{equation}
By Lemma~\ref{th:factorization}, $s_{D\cap B(y,1)}(v)$ is comparable to
$G_D(x_0,v)$ at $v=y$, thus
\eqref{eq:inaccessibility}
is equivalent to 
\begin{equation}
  \label{eq:ginaccessible}
P_D(x_0,y) = \int_\Rd
  G_D(x_0, v) \nu(v, y)dv<\infty\,.  
\end{equation}
Here $y\in \Rd$. If $D$ is non-Greenian, then $s_{D\cap
  B(y,1)}=s_{B(y,1)}$ for every $y\in \Rd$ because $D^c$ is
  polar, and so every point of $D^c$ is accessible from $D$. Thus, for
  general $D$,
\begin{equation}
    \label{eq:obM}
\partial_{M} D \cap \Rd = \set{y \in
  \partial D \, : \, P_D(x_0, y) = \infty}\,.
  \end{equation}
By Fatou's lemma the function $y \mapsto \Lambda_y(G_D(x_0, \cdot))$
is lower semicontinuous. We see that $\partial_M D$ is Borel measurable, and in fact of type
$\mathcal{G}_\delta$.

\proofof{Theorem~\ref{th:martin}} 
Let $D$ be open and Greenian, and let $y \in \partial D$. By
translation invariance, to study $M_D(\cdot,y)$ we may assume with no
loss of generality, that $y=0$.

Let $x\in D$, $\rho=(|x| \wedge |x_0|)/2$, and $D_\rho=D\cap B(0,\rho)$. 
By Harnack inequa\-li\-ty in {the first} variable,
$G_D(x,v)/G_D(x_0,v)$ is bounded from above and below for $v\in D_\rho$.
Also, $G_D(x,v)=P_{D_\rho}[G(x,u)du](v)$,
$G_D(x_0,v)=P_{D_\rho}[G(x_0,u)du](v)$ for $v\in D_\rho$.
Lemma~\ref{th:oscillation}, applied to $D_\rho$, yields that $M_D(x,0)$ 
is well-defined by (\ref{eq:martin:definition}). Clearly, $0<M_D(x,0)<\infty$.

Denote $M(x)=M_D(x, 0)$. If $0$ is inaccessible from $D$, then by \eqref{eq:limit:explicite} we have
$$
  M(x) = 
  \int_D \nu(0, y) G_D(x, y) dy \Big/ \int_D \nu(0, y) G_D(x_0, y) dy =
  P_D(x, 0) / P_D(x_0, 0) \, ,
$$
in particular $M(x)$ {is $\alpha$-harmonic on $D$ with outer charge 
$(P_D(x_0, 0))^{-1} \varepsilon_0$ and so} it is {\it not} 
$\alpha$-harmonic on $D$ {with zero outer charge}. 
However, if $0$ is accessible from $D$ then 
\begin{equation}\label{eq:Mharm}
M(x)=\int_{D\setminus U}M(y)\omega^x_U(dy)\,,\quad x\in U\,,
\end{equation}
for every  $U = D \setminus \cl{B}_R$ with $R>0$. 
Indeed, \eqref{eq:Mharm} is equivalent to uniform integrability of
$G_D(y,z)/G_D(x_0,z)$ with respect to $\omega^x_U(dy)$ on the (bounded)
set $D\setminus U$ as $D\ni z \to 0$. To prove the uniform
integrability, let $0<r<\min(R/4,|x_0|/4)$ and $z_0 \in D_r$ be a fixed point. For
$y \in D_R \setminus D_{3r}$ and $z\in D_r$, Remark~\ref{rem:Greenfunction}
yields that 
$$
   \frac{G_D(y, z)}{G_D(x_0, z)} \leq
   C_{d, \alpha, r} \frac{G_D(y,z_0)}{G_D(x_0,z_0)}.
$$
Again by Remark~\ref{rem:Greenfunction} we obtain that $\sup_{y \in D_R \setminus D_{3r}} G_{D}(y,z_0) < \infty$.
Thus we only need to estimate $\int_{D_{3 r}} G_D(z, y)
\omega^x_U(dy)/G_D(x_0, z)$ for $z\in D_r$.

Since the density function (Poisson kernel) of $\omega^x_U$ is bounded
on $D_{3 R/4}$, we have
\begin{equation} \label{eq:martin:accessible:3}
  \int_{D_{3 r}} G_D(y, z) \omega^x_U(dy) \leq
  C_{d, \alpha, D, R} \, \int_{D_{3 r}} G_D(y, z)dy\,.
\end{equation}
By (\ref{eq:omega:density}), $\omega^{x_0}_{D \setminus \overline{D}_{3 r}}$ is
absolutely continuous on $\overline{D}_{3r}$ with respect to the
Lebesgue measure, and has $P_{D \setminus \cl{D}_{3 r}}(x_0, \cdot)$
as density function. Thus
\begin{eqnarray}
  G_D(x_0, z)
  & = & \nonumber
  \int_{D_{3 r}} G_D(y, z) P_{D \setminus \cl{D}_{3 r}}(x_0, y) dy
  \\ & = & \nonumber
  \int_{D_{3 r}} \int_{D \setminus \cl{D}_{3 r}} G_D(y, z)
    G_{D \setminus \cl{D}_{3 r}}(x_0, \zeta) \nu(\zeta, y) d\zeta dy
  \\ & \geq & \label{eq:martin:accessible:5}
  2^{-d-\alpha}\left( \int_{D_{3 r}} G_D(y, z) dy \right)
    \left( \int_{D \setminus \cl{D}_{3 r}}
    G_{D \setminus \cl{D}_{3 r}}(x_0, \zeta) \nu(\zeta, 0) d\zeta \right) \,.
\end{eqnarray}
The last integral becomes arbitrarily large when $r$ is small
enough. This is because $\int_D G_D(x_0, \zeta) \nu(\zeta,
0)d\zeta=\infty$, $0$ is accessible from $D$, and $G_D(x_0, \cdot)\approx s_{D {\cap B}}(\cdot)$ at $0$ by
Lemma~\ref{th:factorization}. 

Combining this, \eqref{eq:martin:accessible:5}, and
\eqref{eq:martin:accessible:3}, we obtain the uniform integrability, and \eqref{eq:Mharm}.
In fact, (\ref{eq:omega:harmonic}) yields \eqref{eq:Mharm} for every
open $U\subset D$ provided $0 \notin \cl{U}$. In particular, $M$ is
$\alpha$-harmonic on $D$. 
Regarding the remark at the end of Section~\ref{sec:harmonicity} we note that $f=M$ violates
(\ref{eq:rh}) because $M$ vanishes on $D^c$.

We now turn to the Martin kernel with the pole at infinity. Let $x\in D$.
If $D = \Rd$ and $\Rd$ is Greenian, or $\alpha<d$, then 
$M_D(x, \infty) = \lim_{|v|\to \infty}|v-x|^{\alpha-d}/|v-x_0|^{\alpha-d}=1$, $s_D\equiv\infty$, 
and we are done. 
Without loosing generality we may suppose in what follows that $D$ is a proper
unbounded (Greenian) subset of $\Rd$, and $0 \in D^c$. 
Consider the inversion with respect to the unit sphere: 
$$
Tx = \frac{1}{|x|^2}x\,,\quad x\neq 0\,.
$$
Inversion can be used to reduce potential theoretic problems
for the point at infinity to those at $0$, see \cite{bib:bz} for a
detailed discussion.
As proved in \cite{bib:bz}, the set $T D=\{Tx:\;x\in D\}$ has Green function
\begin{equation}\label{eq:tkfg}
G_D(x, v) = |x|^{\alpha-d}|v|^{\alpha-d} G_{T D}(Tx, Tv)\,,\quad x,
v\neq 0\,,
\end{equation}
in particular it is Greenian.
We obtain
\begin{equation} \label{eq:martin:kelvin}
M_D(x,\infty)
=
\lim_{D\ni v\to \infty}\frac{|x|^{\alpha-d}|v|^{\alpha-d}G_{TD}(Tx,Tv)}
{|x_0|^{\alpha-d}|v|^{\alpha-d}G_{TD}(Tx_0,Tv)}
=\frac{|x|^{\alpha-d}}{|x_0|^{\alpha-d}}M_{TD}(Tx,0)\,.
\end{equation}
Here $M_{TD}$ denotes the Martin kernel of $TD$ with the pole at $0$ and the reference point at $T
x_0$. The existence of $M_D(x,\infty)$ defined by
\eqref{eq:martin:definition} is proved.
Also, $0<M_D(x,\infty)<\infty$.
Note that $|x|^{\alpha-d}M_{TD}(Tx,0)$ in (\ref{eq:martin:kelvin}) is
the {\it Kelvin transform}\/ of $M_{TD}(x,0)$, see \cite{bib:bz}. 
By \cite{bib:bz}, $M_D(x,\infty)$ is $\alpha$-harmonic
in $D$ if and only if $M_{TD}(x,0)$ is $\alpha$-harmonic in $TD$.
We finally observe that $\infty$ is accessible for $D$ if and only if $0$ is accessible for
$TD$. Indeed, by (\ref{eq:tkfg}) and a change of variable $v=Ty$ with
Jacobian $|y|^{-2d}$,
\begin{equation} \label{eq:kelvin:green}
  \int G_{TD}(Tx, y) \nu(0, y) dy =
  \mathcal{A}_{d, -\alpha} |x|^{d-\alpha} \int G_D(x, Ty) |y|^{-2d} dy =
  \mathcal{A}_{d, -\alpha} |x|^{d - \alpha} s_D(x) \,.
\end{equation}
Therefore $\alpha$-harmonicity of $M_D(x,\infty)$ is equivalent to
accessibility of $\infty$ for $D$.
\qed

We let 
\begin{equation}
  \label{eq:dm}
M_D(x, y) = \frac{G_D(x, y)}{G_D(x_0, y)}\,,\quad x,y \in D\,,
\end{equation}
so that $M_D(x, y)$ is now defined for all $x \in D$ and $y \in
D^*$. 
In passing we note that (\ref{eq:tkfg}) yields 
\begin{equation}
  \label{eq:ktMk}
  M_D(x,y)=\frac{|x|^{\alpha-d}}{|x_0|^{\alpha-d}}M_{TD}(Tx,Ty)\,,\quad
  x\in D\,,\quad y\in D^*\,,
\end{equation}
where $0\notin D$, the reference point for $M_{TD}$ is $Tx_0$ and we use the
convention $T0=\infty$, $T\infty=0$.
Recall that $B_r=B(0,r)$.
\begin{lemma} \label{th:martin:continuity}
For every $\rho>0$ and $\eta > 0$ there is $r > 0$ such that for every Greenian $D$
\begin{equation} \label{eq:martin:oscillation}
  \ro_{y\in \cl{D} \cap B_r} M_D(x, y) \leq
  1 + \eta
  \, , \quad \mbox{if } \;\; x, x_0 \in D \setminus \cl{B}_\rho \, ,
\end{equation}
\begin{equation} \label{eq:martin:oscillation:infinity}
  \ro_{y\in D^* \setminus \cl{B}_{1/r}} M_D(x, y) \leq
  1 + \eta
  \, , \quad \mbox{if } \;\; x, x_0 \in D \cap \cl{B}_{1/\rho} \,.
\end{equation}
The Martin kernel $M_D(x, y)$: $D \times D^* \setminus \{(x_0,
x_0)\}\mapsto [0,\infty]$ is jointly continuous.
\end{lemma}
\proof
To prove (\ref{eq:martin:oscillation}) let $\rho>r>0$ and $x, x_0 \in D \setminus \cl{B}_\rho$. We note
that by (\ref{eq:martin:definition})
$$
\sup_{y\in \cl{D}\cap {B}_r}M_D(x,y)=
\sup_{y\in D\cap {B}_r}\frac{G_D(x,y)}{G_D(x_0,y)}\,,\quad
\inf_{y\in \cl{D}\cap {B}_r}M_D(x,y)=
\inf_{y\in D\cap {B}_r}\frac{G_D(x,y)}{G_D(x_0,y)}\,,
$$
hence
$ \ro_{y\in \cl{D}\cap {B}_{r}} M_D(x, y)=
 \ro_{y\in D \cap {B}_{r}} M_D(x, y)$.
As functions of $y$, $G_D(x, y)$ and $G_D(x_0, y)$ are nonnegative Poisson
integrals on $D\cap B_\rho$ of measures on $B_\rho^c$. Thus
(\ref{eq:martin:oscillation}) is an immediate consequence of
Lemma~\ref{th:oscillation} and scaling. 
To prove \eqref{eq:martin:oscillation:infinity}, as in the proof of
Theorem~\ref{th:martin} we may assume 
that $0 \in D^c$, and then \eqref{eq:martin:kelvin} reduces 
\eqref{eq:martin:oscillation:infinity} to
\eqref{eq:martin:oscillation} for $TD$.

By Lemma~\ref{lem:cgf}, $M_D$ given by \eqref{eq:dm} is jointly continuous: 
$D\times D\setminus \{(x_0,x_0)\}\mapsto [0,\infty]$.
We will consider the remaining case; let $D\times D^*\ni (x',y')\to (x,y) \in D\times
\partial_* D$. We have
$$
  \frac{M_D(x', y')}{M_D(x, y)} =
  \frac{M_D(x', y')}{M_D(x, y')} \cdot \frac{M_D(x, y')}{M_D(x, y)}\,.
$$
Here the second factor on the right hand side converges to $1$ by
\eqref{eq:martin:oscillation} or
\eqref{eq:martin:oscillation:infinity}.
We will verify uniform continuity of the first factor at $x'=x$.
If $\cl{B}(x,s)\subset D$, and $y'\in \cl{B}(x,s)^c$ then by (\ref{eq:green:harmonic}),
\eqref{eq:Mharm} and \eqref{eq:harm} we see that
{$f(\cdot)=M_D(\cdot, y')$}  
satisfies \eqref{eq:harmonic:definition} {with $\lambda(dw) =
  M_D(w,y')dw+c\varepsilon_{y'}(dw)$ and $U$ being the open ball}.  Here
$c=1/ P_D(x_0, y')$ if $y'\in D^c$ and $y'$ is inaccessible for $D$,  and $c=0$ otherwise.
The uniform continuity follows from \eqref{eq:cm} as in the proof of
Lemma~\ref{lem:cgf}, and $M_D(x', y') /M_D(x, y) \rightarrow 1$.
\qed
\begin{remark}\label{rem:ha}
{\rm
The kernel functions $G_D$, $P_D$, and $M_D$ may be studied
without explicit mention of harmonicity, by using
(\ref{eq:green:twosets}) and its consequences, (\ref{eq:harm}) and (\ref{eq:Mharm}).
}
\end{remark}


\section{Structure of nonnegative harmonic functions} \label{sec:harmonic}

\begin{lemma} \label{th:martin:decomposition}
If $f\geq 0$ and $f$ is $\alpha$-harmonic on $D$, with outer charge $\lambda\geq 0$, then
there is a unique {nonnegative} function $f_s$ $\alpha$-harmonic in
$D$ {with zero outer charge} such that 
$f = P_D[\lambda]+f_s$ on $D$.
\end{lemma}
\proof
Let $D_n$ be an increasing sequence of open precompact subsets of $D$
such that $\bigcup_{n = 1}^{\infty} D_n = D$.  
By \eqref{eq:poisson:definition}, monotone convergence of $G_{D_n}$ to
$G_D$, $\alpha$-harmonicity of $f$, and \eqref{eq:n}, we have
$$
  {P}_D[\lambda](x) =
\lim_{n \rightarrow \infty} \int_{D^c} \int_D 
G_{D_n}(x, z) \nu(z, y) dz {\lambda}(dy) \leq f(x) \,,\quad x\in D\,. 
$$
We let $f_s = f - P_D[\lambda]$. 
The stated properties easily follow. 
\qed
%

Since the outer charge of $f_s$ vanishes on $D^c$, the present setting for further decomposition of $f_s$ is
analogous to those of \cite{{bib:b:repr},{bib:ms}, bib:sw}, despite the initial generality of our definition
of harmonic functions (but see the discussion following (\ref{eq:representation})). 

\begin{lemma} \label{lem:martin:representation}
Let $D$ be Greenian and let $\mu \geq 0$ be a finite measure on $\partial_M D$. Then
\begin{equation} \label{eq:martin:representation}
  f(x) =
  \int_{\partial_M D} M_D(x, y) \mu(dy)
  \, , \quad x \in \Rd
\end{equation}
is $\alpha$-harmonic on $D$ {with zero outer charge}. Conversely, if
 $f\geq 0$ is $\alpha$-harmonic on $D$ {with zero outer charge} then there
is a unique finite measure $\mu\geq 0$ on $\partial_M D$ satisfying 
\eqref{eq:martin:representation}.
\end{lemma}
\proof
It is a straightforward consequence of Theorem~\ref{th:martin} that
$f$ given by \eqref{eq:martin:representation} is $\alpha$-harmonic in
$D$ with zero outer charge on $D^c$.
We shall write $f=M_D[\mu]$.

Let $f$ be a nonnegative function $\alpha$-harmonic in $D$ with zero
outer charge on $D^c$. Let $D_n$ denote an increasing sequence of open sets
precompact in $D$ such that $\bigcup_{n = 1}^{\infty} D_n = D$.
We will also assume that $\omega^x_{D_n}(\partial
D_n) = 0$ for $x \in D_n$, which holds, e.g., if $D_n$ are
Lipschitz domains. By \eqref{eq:omega:density} for $x \in D_n$ we then have 
$$
  f(x) =
  \int_{D \setminus D_n} P_{D_n}(x, y) f(y) dy =
  \int_{D_n} M_{D_n}(x, v) \left( G_{D_n}(x_0, v) \int_{D \setminus D_n} \nu(v, y) f(y) dy \right) dv
  \,. 
$$
Let $\mu_n(dv) = \left(G_{D_n}(x_0, v) \int_{D \setminus D_n} \nu(v,
  y) f(y) dy \right) dv$. Since $\mu_n(D) = f(x_0)<\infty$, by considering a subsequence
we may assume that $\mu_n$ weakly converge to a finite nonnegative measure $\mu$ on
$D^*$. We claim that $\mu$ is supported in $\partial_* D$.
Indeed, if $n > k$, $v \in D_k$, $y \in D \setminus D_n$, then
$\nu(v, y) \leq C_k$ and $G_{D_n}(x_0, v) \leq G_D(x_0, v)$. Hence
$$
  \mu_n(D_k) \le
  C_k \left( \int_{D_k} G_D(x_0, v) dv \right)
    \left( \int_{D \setminus D_n} f(y) (1+|y|)^{-d-\alpha}dy \right) \rightarrow 0
$$
as $n \rightarrow \infty$, see \eqref{eq:calkfh}. 
This proves that $\mu(D_k) = 0$ and so $\mu$ is a measure on $\partial_* D$.

Let $\varepsilon > 0$ and $x \in D$. 
By Lemma~\ref{th:martin:continuity} for every $y \in \partial_* D$ there exists its
neighborhood $V_y$ such that
$$
  \ro_{V_y} M_U(x, \cdot) \leq
  1 + \varepsilon\,,
$$
with $U=D$ {\it and} $U=D_n$, $n=1,\ldots$.
From these, one selects a finite family
$\{V_j,\,j = 1, \ldots, m\}$ such that $V = V_1 \cup \ldots \cup V_m \supset \partial_* D$.
For $j=1,\ldots,m$, let $z_j \in V_j \cap D$.
Let $k$ be so large that for $n \geq k$ we have
$z_j \in D_n$, and 
$$
  (1 + \varepsilon)^{-1} \leq
  \frac{M_D(x, z_j)}{M_{D_n}(x, z_j)} \le 1 + \varepsilon
  \, , \quad j = 1, \ldots, m \,.
$$
If $v \in V_j \cap D_n$ 
then
$$
  (1 + \varepsilon)^{-3} \leq
  \frac{M_D(x,v)}{M_D(x,z_j)} \cdot
    \frac{M_D(x,z_j)}{M_{D_n}(x,z_j)} \cdot
    \frac{M_{D_n}(x,z_j)}{M_{D_n}(x,v)} \leq
  (1 + \varepsilon)^3 \,.
$$
Therefore
$$
  (1 + \varepsilon)^{-3} \leq
  \frac{\int_{D \cap V} M_D(x, y) \mu_n(dy)}
    {\int_{D \cap V} M_{D_n}(x, y) \mu_n(dy)} \leq
  (1 + \varepsilon)^3
  \, , \quad n \geq k \,.
$$
By letting $n \rightarrow \infty$ we obtain
$$
  (1 + \varepsilon)^{-3} \leq
  \frac{\int_{\partial_* D} M_D(x, y) \mu(dy)}{f(x)} \leq
  (1 + \varepsilon)^3 \, ,
$$
which gives $f(x)=\int_{\partial_* D}M_D(x,y)\mu(dy)$.

We will prove that $\mu$ is concentrated on $\partial_M D$. Let $U$ be
open and precompact in $D$ and let $x \in U$.  By Theorem~\ref{th:martin}
and \eqref{eq:harm}, 
if $y \in \partial_* D$, then $M_D(x, y) \geq \int_{D \setminus U}
M_D(z, y) \omega^x_U(dz)$ and equality holds if and only if $y \in
\partial_M D$. By Fubini's theorem
$$
  0 =
  f(x) - \int_{D \setminus U} f(z) \omega^x_U(dz) =
  \int_{\partial D} \expr{M_D(x, y) - \int_{D \setminus U} M_D(z, y) \omega^x_U(dz)} \mu(dy) \, ,
$$
hence $\mu(\partial_* D \setminus \partial_M D) = 0$.

We will prove the uniqueness of $\mu$ in the representation \eqref{eq:martin:representation}.
We first consider $f =M_D[\varepsilon_{y_0}]= M_D(\cdot, y_0)$, where $y_0 \in \partial_M D$. 
To simplify notation, we assume as we may that $y_0 = 0$ 
(we use translation invariance if $0\neq y_0\in \Rd$ and inversion if $y_0 = \infty$). 

Let $D_r = D \cap B_r$, $D_r' = D \setminus \cl{B}_r$.
Suppose that $f$ satisfies \eqref{eq:martin:representation} for a
nonnegative measure $\mu$ on $\partial_M D$. Let $r > 0$ and 
$g(x) =\int_{|y|>3r} M_D(x, y) {\mu}(dy)$. 
Considering $y\in \partial_M D$ such that $|y|>3r$, by \eqref{eq:Mharm} we get
$$
  g(x) =
  \int_{D \setminus D_{2 r}} g(z) \omega^x_{D_{2 r}}(dz)
  \, , \quad x \in D_{2 r}\,.
$$
On the other hand, we may apply Lemma~\ref{th:omega:twofunctions} to
$f$, $g$, and $D_r'$, to verify that
$$
  g(x) =
  \int_{D \setminus D_r'} g(z) \omega^x_{D_r'}(dz)
  \, , \quad x \in D_r'\,.
$$
Lemma~\ref{th:omega:twosets} yields $g(x) = \int_{D^c}
g(z) \omega^x_D(dz) = 0$, that is, $\mu = 0$ on $\partial_M D\cap
\{|y|>3r\}$. 
In particular, the measures $\mu_n$ considered at the
  beginning of the proof, corresponding to
$f(\cdot)=M_D(\cdot,y_0)$, weakly converge to $\varepsilon_{y_0}$.
Fubini's theorem and dominated convergence yield that for general $f=M_D[\mu]$
the measures $\mu_n$ corresponding to $f$ weakly converge to
$\mu$. Since $\mu_n$ are determined by $f$, so is $\mu$.
\qed

We note that if $f$ is $\alpha$-harmonic in $D$ {(with zero outer charge)} and $0 \leq f \leq
M_D(\cdot, y_0)$ then the proof of Lemma~\ref{lem:martin:representation} yields  
$f = c\, M_D(\cdot, y_0)$ for some $c \in [0, 1]$. 
Thus, $M_D(\cdot,y_0)$ is {\it minimal harmonic}\/ i.e.\ an extremal
point of the class of nonnegative functions $f$
$\alpha$-harmonic on $D$ {(with zero outer charge)}, such that
$f(x_0)=1$. 
We note, however, that our proof 
of Lemma~\ref{lem:martin:representation}
does not invoke Choquet's
theorem. Instead it relies on (\ref{eq:martin:definition})
and Lemma~\ref{th:oscillation}.

\proofof{Theorem~\ref{th:representation}}
The theorem collects results of 
Lemma~\ref{th:martin:decomposition}~and~\ref{lem:martin:representation}.
\qed

In particular, if the point at infinity is inaccessible for $D$ then
$\infty$ is not charged by the measure $\mu$ in the representation 
\eqref{eq:martin:representation}, and $M(\cdot,\infty)=s_D$ is not
$\alpha$-harmonic in $D$ with zero outer charge, compare (\ref{eq:superharmonic}).


\section{Miscelanea} \label{sec:miscelanea}


Consider $f(x) = \omega^x_D(\partial_M D)$, $x \in \Rd$. By Lemma~\ref{th:martin:decomposition},
$$
  \int_{\partial_M D} P_D(x, y) dy \leq f(x)\leq 1\, , \quad x \in D \,.
$$
Since $P_D(x, y) = \infty$ for $y \in \partial_M D$, we conclude that $|\partial_M D| = 0$.

We will now strengthen the result of Lemma~\ref{th:poisson:harmonic} and \eqref{eq:omega:density}.
\begin{proposition} \label{th:omega:poisson}
For Greenian $D \subset \Rd$, and $x\in D$, the harmonic measure $\omega^x_D$ 
is absolutely continuous on $D^c \setminus \partial_M D$ with respect
to the Lebesgue measure, with density $P_D(x, \cdot)$.
\end{proposition}
\proof
Let $K \subset D^c \setminus \partial_M D$ be compact and let $f(x) =
\omega^x_D(K) - P_D[\ind_K](x)\geq 0$. We will verify that $f = 0$.
By Theorem~\ref{th:representation}, $f(x) = \int_{\partial_M D}  M_D(x,
y) \mu(dy)$ for some nonnegative finite measure $\mu$ on $\partial_M D$. Let
$L \subset \partial_M D$ be compact and let $g(x) = \int_L M_D(x, y)
\mu(dx)$. It suffices to prove that $g = 0$. We let
$$
  U =
  \{x \in D \, : \, 2 \dist(x, K) \leq \dist(x, L) \}
  \, , \quad
  V =
  \{x \in D \, : \, 2 \dist(x, L) \leq \dist(x, K) \} \,.
$$
Observe that by \eqref{eq:Mharm}, $g(x) = \int_{D \setminus U} g(y)
\omega^x_U(dy)$ for $x \in U$. On the other hand,
Lemma~\ref{th:omega:twofunctions} applied to $\omega^x_D(K)$, $g$, and
$V \subset D$ yields $g(x) = \int_{D \setminus V} g(y)
\omega^x_V(dy)$ for $x \in V$. Hence we may apply
Lemma~\ref{th:omega:twosets} to conclude that $g(x) = 0$.
\qed

In particular, if $f = P_D[\lambda]$, where $\lambda$ is nonnegative 
and absolutely continuous with respect to the Lebesgue measure 
on $D^c$ and has a density function $g$, then we can write 
$$
  f(x) =
  \int_{D^c} g(y) \omega^x_D(dy)\,,\quad x\in D \,.
$$
This, however, requires a {\it convention}\/ that $g(y) = 0$ for $y \in
\partial_M D$ on the right hand side, and should be used with
caution. Another common convention is writing $f$ instead of $g$
above, see (\ref{eq:harmonic:definition-ac}).

We note that there are domains $D$ for which the part of the harmonic
measure which is singular with respect to the Lebesgue measure
(i.e.  $\omega^x_D$ on $\partial_M D$) is positive. 
Indeed, such is the complement of every closed non-polar set of zero
Lebesgue measure, for example, the complement of a point on the line if 
$1 <\alpha < 2$, see \cite{bib:po}.

\begin{lemma}\label{l:hc}
Every nonnegative $f$ harmonic on non-Greenian $D$
is constant on $D$.
\end{lemma}
\proof
If $\alpha<d$ then $G_D$ is majorized by the Riesz kernel 
\cite{bib:La} and so every (open) $D\subset \Rd$ is Greenian. For $\alpha\geq d=1$, by \cite{bib:po}, if $D$ is
non-Greenian then $D^c$ is polar. In this case, let $x, y\in D$ and $0<r<\min({\rm
  dist}(y,D^c), |y-x|)$. By recurrence (see \cite{bib:po} for the
definition) for every $\varepsilon>0$ there is an open precompact $B\subset D$ such that $x\in
B$ and
$\omega^x_{B\setminus\overline{B}(y,r)}(B(y,r))>1-\varepsilon$. 
Using small $\varepsilon$ and $r$, and continuity of $f$ at $y$ we obtain $f(x)\geq f(y)$, hence $f$ is
constant on $D$.
\qed

We will give examples of accessible and inaccessible boundary points.
Let $d \geq 2$ and let $f : (0, 1) \rightarrow (0, \infty)$ be any
bounded increasing function. 
We define a {\it thorn} $D_f$ by (cf. \cite{bib:bk}):
$$
  D_f =
  \{(x_1, \dots, x_d) \in \Rd \, : \, 0 < x_1 < 1, \, |(x_2, \dots, x_d)| < f(x_1)\} \,.
$$
\begin{proposition} \label{th:thorn}
The origin is inaccessible from $D_f$ if and only if 
$\int_0^1 t^{-d - \alpha} f(t)^{d + \alpha - 1} dt<\infty$.
\end{proposition}
\proof
We denote the above integral by $I_f$.
We need to prove that $\Lambda_0(s_{D_f})=\infty$ if and only if $I_f=\infty$.
Let $g(t) = \frac{1}{2} (f(t / 2) \wedge t)$. 
Note that $I_f = \infty$ if and only if $I_g = \infty$ 
(see the proof of Theorem~1.1(i) in \cite{bib:bk}). 
For small $x \in D_g$ we have $B(x, g(x_1)) \subset D_f$,
hence
$$
  s_{D_f}(x) \geq
  s_{B(x, g(x_1))}(x) =
  C_{d, \alpha} (g(x_1))^\alpha.
$$
Thus, if $I_f=\infty$ then ($I_g=\infty$ and) $\Lambda_0(s_{D_f}) = \infty$.
We may now assume that $I_f$ is finite and $f(t) \leq |t|$. 
Let $D_{f, r} = D_f \cap B_r$ for $r>0$. Let $r<1/4$. We have
$$
  s_{D_f}(x) =
  s_{D_{f, 4 r}}(x) +
    \int_{D_f \setminus D_{f, 4 r}} s_{D_f}(y) \omega^x_{D_{f, 4
  r}}(dy)\,.
$$
The latter term is a Poisson integral on $D_{f, 4 r}$. In view of Lemma~\ref{th:factorization}
$$
  s_{D_f}(x) \leq
  s_{D_{f, 4 r}}(x) (1 + C_{d, \alpha} \Lambda_{0, 3 r}(s_{D_f}))
  \, , \quad x \in D_{f, 3 r} \,.
$$
Let $M(r) = \sup\limits_{x_1 = r} s_{D_f}(x) / (f(4 r))^\alpha$. 
Inscribing $D_{f, r}$ into a cylinder and observing that 
$G_D(x, y) \leq C_{d, \alpha} |x - y|^{-d + \alpha}$ 
we can show that $s_{D_{f, r}}(x) \leq C_{d, \alpha}
(f(r))^\alpha$. We thus obtain
$$
  M(r) \leq
  c_1 + c_2 \int_{2 r}^1 M(t) (f(4 t))^{d + \alpha - 1} t^{-d - \alpha} dt \, ,
$$
where $c_1$ and $c_2$ are some constants depending on $d$ and $\alpha$. Let $R > 0$ satisfy
$$
  2 c_2 \int_0^R (f(t))^{d + \alpha - 1} t^{-d - \alpha} dt <
  1 \,,
$$
so that
$$
  M(r) \leq
  c_1 + (1/2) \sup_{(2 r, R)} M + c_2 I_f \sup_{(R, 1)} M \,.
$$
It follows that $M$ is bounded by $2 c_1 + 2 c_2 I_f \sup_{(R,
  1)} M$. By using the definition of $M$, we conclude that $\Lambda_0(s_{D_f})$ is finite.
\qed

We note that by Fatou's lemma, if $y \in \partial_M D$
then $P_D(x, z) \rightarrow P_D(x, y) = \infty$ as $\cl{D}^c \ni z
\rightarrow y$. If $y \in \partial D\setminus \partial_M D$ 
and $\cl{D}^c \ni z \rightarrow y$ non-tangentially 
(i.e. $|z - y| \leq c \dist(z,D)$ for some $c > 0$)
then by dominated convergence we have $P_D(x, z) \rightarrow P_D(x, y)
< \infty$. 

The next result is an extension of \cite[Lemma~7]{bib:b:repr}.
\begin{proposition}\label{l:jpjm}
If $y \in \partial_M D\cap ((\cl{D})^c)^*$ then 
\begin{equation}
  \label{eq:wjm}
  M_D(x, y) =
  \lim_{\cl{D}^c \ni z \rightarrow y} \frac{P_D(x, z)}{P_D(x_0, z)}
  \,.
\end{equation}
\end{proposition}
\proof
Suppose that $y=0$ is a limit point of $D$ and
of the interior of $D^c$, and $0$ is accessible for $D$. Let $x\in D$.
Recall our notation $D_r = D \cap B_r$, $D_r' = D \setminus D_r$. 
Let $0<4 r < |x| \wedge |x_0|$, and let $z \in B_r\setminus \cl{D}$. 
By Lemma~\ref{th:factorization} and Remark~\ref{rem:sf}
$$
  \int_{D_r} G_D(x, v) \nu(v, z) dv \geq
  C_{d, \alpha} \int_{D_{3 r}} s_{D_{2 r}}(v) \nu(v, y) dv 
    \int_{D \setminus D_{2 r}} G_D(x, v) \nu(v, y) dv \,.
$$
This also holds for $x=x_0$. By Fatou's lemma we have 
$\lim_{\cl{D}^c\ni z \rightarrow 0} \int_{D_r} s_{D_{2 r}}(v) \nu(v, z) dv= \infty$.
Since $\int_{D \setminus D_{r}} G_D(x, v) \nu(v, z) dv$ is bounded in
$z$, \eqref{eq:poisson:definition} yields
$$
  \lim_{\cl{D}^c \ni z \rightarrow 0} \frac{P_D(x, z)}{P_D(x_0, z)} =
  \lim_{\cl{D}^c \ni z \rightarrow 0}
    \frac{\int_{D_r} G_D(x, v) \nu(v, z) dv}{\int_{D_r} G_D(x_0, v) \nu(v, z) dv} \, ,
$$
provided that limits exist. If $\delta > 0$ then for sufficiently
small $r$ by \eqref{eq:martin:definition} we obtain
$$
  M_D(x, 0) - \delta \leq
  \frac{\int_{D_r} G_D(x, v) \nu(v, z) dv}{\int_{D_r} G_D(x_0, v) \nu(v, z) dv} \leq
  M_D(x, 0) + \delta \, ,
$$
which proves \eqref{eq:wjm}.
For general $y\in \partial D$ we use translation invariance. 
If $y=\infty$ then we use inversion.
Namely, (\ref{eq:poisson:definition}) and \eqref{eq:tkfg} and
$|Tx-Tz|=|x-z|/(|x||z|)$ lead to 
$$P_D(x, z) =|x|^{\alpha - d} |z|^{-\alpha - d} P_{TD}(Tx, Tz)\,,$$
see \cite{bib:bz}. This, and \eqref{eq:martin:kelvin} yield \eqref{eq:wjm}. 
\qed

If $D=B(0,r)$, $r>0$, and $x_{0}=0$, then we have
\begin{equation}\label{eq:wjmk}
M_D(x,y)=
r^{d-\alpha}
\frac{\left(r^{2}-|x|^{2}\right)^{\alpha/2}}{|x-y|^{d}}\,,\quad |x|<r\,,
\end{equation}
for every $y\in \partial B(0,r)$. (\ref{eq:wjmk}) follows from
Proposition~\ref{l:jpjm} and (\ref{eq:poisson:ball}) or (\ref{eq:martin:definition}) and (\ref{wfg}). 
The formula 
was given before in \cite{bib:hm}, \cite{bib:b:repr}, \cite{bib:cs}. 
We note that $B_r$ has
all its boundary points $y$ 
accessible
because
$G_{B_r}(x,v)\approx (r-|v|)^{\alpha/2}$ 
as $B_r\ni v\to y$, see  (\ref{wfg}). More generally, a Lipschitz (or even $\kappa$-fat)
domain has all its boundary points
accessible, as follows from
\cite{bib:b:bhp} (\cite{bib:sw}). For more information on the boundary
potential theory in Lipschitz domains we refer to the papers \cite{bib:b:repr},
\cite{bib:babo}, \cite{bib:kmmr}, \cite{bib:sw}, which may suggest further
applications.

\begin{remark}\label{rem:u}
{\rm
The statement of uniqueness in Theorem~\ref{th:representation} may be
strengthened. Namely, nonnegative measures $\mu$ on $\partial_M D$
and $\lambda$ on $D^c_M$ are uniquely determined by the values of
\begin{equation} 
  \int_{D^c_M} P_D(x, y) \lambda(dy) + \int_{\partial_M D} M_D(x, y)
  \mu(dy)\,,\quad \mbox{for }x\in D\,, 
\end{equation}
provided they are finite for (some, hence for all) $x\in D$.
This follows from Theorem~4.2 in \cite{bib:cs}, which states that
a (genuine) function $\alpha$-harmonic on open $U\neq
\emptyset$ is determined $a.e.$ on $\Rd$ by its values on $U$. By
a convolution with smooth compactly supported approximate identity
(integrability follows from (\ref{eq:calkfh})),
this yields uniqueness of $\lambda$ and $\mu$ on $\Rd$. If $\infty\in
\partial_M D$, and $\mu$ has an atom at $\infty$ then the mass of the 
atom is determined by the values of $M_D(x,\infty)\mu(\{\infty\})$.
}
\end{remark}

\begin{remark}\label{rem:gh}
{\rm
Theorem~\ref{th:bhp} extends to general (nonnegative) harmonic
functions such that measures $\lambda$ and $\mu$ in the
representation (\ref{eq:representation}) do not charge $G\cap D^c$ in
Theorem~\ref{th:bhp}.
To be specific, Theorem~\ref{th:bhp} applies to 
$P_D[\lambda]+M_D[\mu]$ in (\ref{eq:representation}).  
Indeed, by (\ref{eq:Mharm}), for finite $y\in \partial_M D$ 
we can consider $M_D(\cdot,y)$ a Poisson integral on a wide class of subdomains 
$U\subset D$ such that $\dist(y,U)>0$, see (\ref{eq:omega:density}) in
this connection. A similar property holds for
$M_D(\cdot,\infty)$ on on a wide class of bounded domains, if $\infty\in \partial_M D$.
This follows from the transformation rules for the harmonic measure
under the inversion of the domain, see \cite{bib:bz}. 

In domain with regular geometry, e.g. Lipschitz, the condition
$(\lambda+\mu)(G\cap D^c)=0$ is equivalent to the traditional
assumption of continuity and vanishing
of a harmonic function on this set, see \cite{bib:b:bhp},
\cite{bib:b:repr} and the references there. Continuous decay of harmonic functions cannot,
however, be required at irregular points of a domain.  
}  
\end{remark}

We finally wish to provide probabilistic interpretations of our results. 
For a general perspective on probabilistic potential theory and for probabilistic
notions mentioned below we refer the reader to, e.g.,
\cite{ChZ}, \cite{bib:bh}, \cite{BBpms2000} and \cite{BB1}.
Here we will only indicate a few specific interpretations as they may 
suggest further extensions. 

The second term on the right hand side of (\ref{eq:representation}) is
coined singular $\alpha$-harmonic  in \cite{bib:b:repr}, \cite{bib:cs} and
\cite{bib:ms} (when set to zero on $D^c$).
As explained in \cite{bib:cs}, in the case of Lipschitz $D$ the function is harmonic for the 
isotropic $\alpha$-stable L\'evy process killed on leaving $D$
(see also \cite{bib:ms}). 
For a general domain $D$ it is more appropriate to relate such functions to 
{\it the continuous exit}\/ of the trajectory of the process from
$D$. 
The observation is implicit in \cite{bib:ms}. 
For example, if $D=B\setminus F$, where $F$
is a non-polar set of Lebesgue measure $0$, then the trajectory of the
process is almost surely continuous when entering $F$. Correspondingly, the harmonic measure
of $F$ with respect to $D$, $x\mapsto \omega_D^x(F)$, is
represented with $\lambda=0$ in (\ref{eq:representation}).

The first term in (\ref{eq:representation}) 
is related to the effect of leaving the domain {\it by a jump}.
The observation is implicit in \cite[formula (5)]{bib:ms}. 
This explains the
role of (\ref{eq:cm}) in our development: the second
term on the right hand side of (\ref{eq:cm}) is the integral
against this part of the harmonic measure, say $\omega^x_{D-}$, which
results from the {\it jumps}\/ of
the trajectory from $D$ to $D^c$, and the Poisson kernel is the
density function of the measure. The latter claim may be verified by
using quasi-left continuity of the process (\cite{bib:b:bhp}).
The reader may want to consider
domains $D$ with $\partial D$ of positive Lebesgue measure, to
apprehend the complexity of the relation between 
$\omega_D$ and $\omega_{D-}$. 
The relation is addressed in Proposition~\ref{th:omega:poisson} above.
In this connection we also refer the reader to \cite{bib:wh} for a
discussion of the important problem of characterization
of domains $D$ for which $\omega_D=\omega_{D-}$.

The above complex behavior (jumps {\it and}\/ continuous exit) 
is manifested only for jump processes, which are exemplified here by the isotropic $\alpha$-stable
L\'evy process on $\Rd$.
In the presence of jumps the distribution of the position of the process stopped when leaving
the domain, i.e.\ the harmonic measure, is supported on $D^c$, but
usually not on $\partial D$, and so it is different from the
distribution of the position of the process immediately before leaving the domain. 
The formula of Ikeda and Watanabe in its full form (for which see, e.g.,
\cite{bib:iw, BB1,BBpms2000}) gives the joint distribution of these two random
variables in terms of the Green function and the L\'evy measure.
Thanks to the simplicity of the L\'evy measure $\nu(x,y)$ in (\ref{eq:n}), 
the estimates for nonnegative $\alpha$-harmonic functions can be
effectively reduced to the estimates of the Green function.
We conjecture that (\ref{eq:representation}) generalizes to a wide class of Markov processes of jump type.

For a Riesz type representation of {\it superharmonic}\/ functions of the
fractional Laplacian on Lipschitz domains we also refer the reader to \cite{bib:cs}.

We {finally} wish to
provide the following probabilistic connection.
The (accessibility) condition
$\Lambda_x(s_D) = \infty$ has appeared implicitly in \cite{bib:bk} and explicitly in
\cite{bib:w}. Authors of these papers consider the following property
of 
{our} symmetric $\alpha$-stable
{L\'evy} process $\{X_t\}$ in $\Rd$ and a given domain $D$: 
{\it There exist{s} a random time interval $(\tau_0, \tau_0 + 1)$
such that $X(t) - X(\tau_0) \in D$ for {all} $t \in (\tau_0, \tau_0 + 1)$.}\/
If $D$ is a {\it thorn}\/ then the property holds if and only if
$\Lambda_0(s_D) = \infty$ (\cite{bib:bk}). 
In \cite{bib:w} all open sets $D$ are considered and the
existence of such interval is established if $\Lambda_0(s_D)$ is
infinite. We conjecture that the {\it accessibility}\/ of $0$ from $D$ is actually a
characterization of this property related to the {\it continuous}\/ 
convergence to $0$ of the trajectories of the corresponding {\it
  conditional}\/ process at its lifetime, see
\cite{BB1}.



{\noindent
{\bf Acknowledgements.} We thank Krzysztof Samotij for his
contribution to our understanding of Martin representation of $\alpha$-harmonic functions.
The first named author gratefully acknowledges the hospitality of the Department of Statistics at Purdue
University, where the paper was written in part. We thank Krzysztof
Burdzy for remarks on accessibility. We thank the referees for many valuable
comments and suggestions.
}


Krzysztof Bogdan (bogdan@pwr.wroc.pl)\\
\noindent
Department of Statistics, Purdue University 
and

\noindent
Institute of Mathematics and Computer Science,

\noindent
Wroc{\l}aw University of Technology    

\noindent
ul. Wybrze{\.z}e Wyspia{\'n}skiego 27, 50-370, Wroc{\l}aw, Poland

Tadeusz Kulczycki (Tadeusz.Kulczycki@pwr.wroc.pl),

Mateusz Kwa\'snicki (Mateusz.Kwasnicki@pwr.wroc.pl)

\noindent
Institute of Mathematics and Computer Science

\noindent
Wroc{\l}aw University of Technology    

\noindent
ul. Wybrze{\.z}e Wyspia{\'n}skiego 27, 50-370, Wroc{\l}aw, Poland


\begin{thebibliography}{[00]}

\bibitem{bib:A}
  H.~Aikawa, 
  \emph{Boundary Harnack principle and Martin boundary for a uniform domain}.
  J. Math. Soc. Japan  53  (2001), no. 1, 119--145.

\bibitem{bib:A2}
  H.~Aikawa,
  \emph{Potential-theoretic characterizations of nonsmooth domains}. 
  Bull. London Math. Soc.  36  (2004), no. 4, 469--482.

\bibitem{bib:AG}
  D.~Armitage, S.~Gardiner, 
  \emph{Classical potential theory}.
  Springer-Verlag, London, 2001.

\bibitem{bib:babo}
  R.~Ba\~nuelos, K.~Bogdan,
  \emph{Symmetric stable processes in cones}. 
  Potential Anal.  21  (2004), no. 3, 263--288.

\bibitem{bib:Ba}
  R.~Bass, 
  \emph{Probabilistic techniques in analysis} 
  Springer-Verlag, New York, 1995.

\bibitem{bib:BB}
  R.~Bass, K.~Burdzy,
  \emph{A probabilistic proof of the boundary Harnack principle}.  
  Seminar on Stochastic Processes, 1989 (San Diego, CA, 1989)
  1--16, Progr. Probab., 18, Birkhäuser Boston, Boston, MA, 1990.

\bibitem{bib:bh}
  J. Bliedtner, W. Hansen,
  \emph{Potential Theory. An analytic and probabilistic approach to balayage}
  Springer-Verlag Berlin Heidelberg, 1986.

\bibitem{BGR}
  R. M. Blumenthal, R. K. Getoor, D. B. Ray.
  \emph{On the distribution of first hits for the symmetric stable processes}.
  Trans. Amer. Math. Soc. 99 (1961), 540--554.

\bibitem{bib:b:bhp}
  K. Bogdan,
  \emph{The boundary Harnack principle for the fractional Laplacian}.
  Studia Math. 123 (1997), 43--80.

\bibitem{bib:b:repr}
  K. Bogdan,
  \emph{Representation of $\alpha$-harmonic functions in Lipschitz domains}.
  Hiroshima Math. J. 29 (1999), 227--243.

\bibitem{Bjmaa2000}
  K.~Bogdan,
  \emph{Sharp estimates for the Green function in Lipschitz domains}.
  J. Math. Anal. Appl. 243 (2000), no. 2, 326--337.


\bibitem{bib:BBC}
  K.~Bogdan, K.~Burdzy, Z.-Q. Chen,
  \emph{Censored stable processes}. 
  Probab. Theory Related Fields  127  (2003), no. 1, 89--152.

\bibitem{BB1} 
  K.~Bogdan and T.~Byczkowski, 
  \emph{Potential theory for the $\alpha$-stable Schr{\"o}dinger operator on bounded Lipschitz domains}. 
  Studia Math{.}, 133 (1999), no. 1, 53--92.

\bibitem{BBpms2000}
  K.~Bogdan, T.~Byczkowski,
  \emph{Potential theory of {S}chr\"odinger operator based on fractional {L}aplacian}.
  Probab. Math. Statist. 20 (2000), no. 2, Acta Univ. Wratislav. No. 2256, 293--335.

\bibitem{bib:bss}
  K.~Bogdan, A.~St\'os, P.~Sztonyk, 
  \emph{Harnack inequality for stable processes on $d$-sets}. 
  Studia Math.  158  (2003), no. 2, 163--198.

\bibitem{bib:bs}
  K.~Bogdan, P.~Sztonyk, 
  \emph{Estimates of potential kernel and Harnack's inequality for anisotropic fractional Laplacian}.
  Preprint (2005) http://arxiv.org/abs/math.PR/0507579

\bibitem{bib:bz}
  K.~Bogdan, T.~\.Zak,
  \emph{On Kelvin transformation}.
  J. Theor. Prob. 19 (2006), no. 1, 89--120.

\bibitem{Br}
  M.~Brelot, 
  \emph{On topologies and boundaries in potential theory}, 
  Lecture Notes in Mathematics, Springer, Berlin, 1971.

\bibitem{bib:bk}
  K. Burdzy, T. Kulczycki,
  \emph{Stable processes have thorns}.
  Ann. Probab. 31 (2003), 170--194  

\bibitem{bib:CK}
  Z.-Q. Chen, P. Kim,
  \emph{Green function estimate for censored stable processes}. 
  Probab. Theory Related Fields  124  (2002), no. 4, 595--610.

\bibitem{bib:cs}
  Z.-Q. Chen, R. Song,
  \emph{Martin boundary and integral representation for harmonic functions
  of symmetric stable processes}.
  J. Funct. Anal. 159 (1998), 267--294.

\bibitem{CS:FK}
  Z.-Q. Chen, R. Song,
  \emph{Conditional gauge theorem for non-local Feynman-Kac transforms}.
  Probab. Theory Relat. Fields 125 (2003), 45--72. 

\bibitem{ChZ}
  K.~Chung, Z.~Zhao, 
  \emph{From Brownian motion to Schr{\"o}dinger's equation}.
  Springer-Verlag, New York, 1995.

\bibitem{bib:Hbhp}
  W.~Hansen,
  \emph{Uniform boundary Harnack principle and generalized
  triangle property}.
  J. Funct. Anal.  226  (2005), no. 2, 452--484.

\bibitem{bib:Hg}
  W.~Hansen,
  \emph{Global comparison of perturbed Green functions}.
  Math. Ann. 334 (2006), no. 3, 643--678. 

\bibitem{bib:hm}
  F.~Hmissi,
  \emph{Fonctions harmoniques pour les potentiels de Riesz sur la
  boule unité} (French) [\emph{Harmonic functions for Riesz potentials on the
  unit ball}]. Exposition. Math.  12  (1994), no. 3, 281--288. 

\bibitem{bib:hw}
  R. A. Hunt, R. L. Wheeden,
  \emph{Positive harmonic functions on Lipschitz domains}.
  Trans. Amer. Math. Soc. 147 (1970), 507--527.

\bibitem{bib:iw}
  N. Ikeda, S. Watanabe,
  \emph{On some relations between the harmonic measure and the L\'evy
  measure for a certain class of Markov processes}.
  Probab. Theory Related Fields 114 (1962), 207--227.

\bibitem{bib:J}
  N.~Jacob,
  \emph{Pseudo differential operators and Markov processes. Vol. I, II, III}.
  Imperial College Press, London, 2001-2005.

\bibitem{bib:TJe}
  T.~Jakubowski,
  \emph{The estimates for the Green function in Lipschitz domains for the
  symmetric stable processes}. 
  Probab. Math. Statist.  22 (2002), no. 2, Acta Univ. Wratislav. No. 2470, 419--441.

\bibitem{bib:Kg}
  T.~Kulczycki,
  \emph{Properties of Green function of symmetric stable processes}. 
  Probab. Math. Statist.  17  (1997), no. 2, Acta Univ. Wratislav. No. 2029, 339--364.

\bibitem{bib:ki}
  T.~Kulczycki, 
  \emph{Intrinsic ultracontractivity for symmetric stable processes}. 
  Bull. Polish Acad. Sci. Math.  46  (1998), no. 3, 325--334.

\bibitem{bib:KW}
  H.~Kunita, T.~Watanabe,
  \emph{Markov processes and Martin boundaries I}.
  Illinois J. Math.  9  (1965), 485--526.

\bibitem{bib:La}
  N. S. Landkof, 
  \emph{Foundations of modern potential theory}.
  Springer-Verlag, New York-Heidelberg 1972.

\bibitem{bib:kmmr}
  K.~Michalik, M.~Ryznar,
  \emph{Relative Fatou theorem for $\alpha$-harmonic
  functions in Lipschitz domains}.
  Illinois J. Math.  48 (2004), no. 3, 977--998.

\bibitem{bib:ms}
  K. Michalik, K. Samotij,
  \emph{Martin representation for $\alpha$-harmonic functions}.
  Probab. Math. Statist. 20 (2000), 75--91.

\bibitem{bib:pi}
  R.~Pinsky,
  \emph{Positive harmonic functions and diffusion}
  Cambridge Studies in Advanced Mathematics, 45. Cambridge University
  Press, Cambridge, 1995.

\bibitem{bib:po}
  S.~Port,
  \emph{Hitting times and potentials for recurrent stable processes}. 
  J. Analyse Math.  20  (1967), 371--395.

\bibitem{bib:Rm}
  M.~Riesz, 
  \emph{Int{\'e}grales de {R}iemann-{L}iouville et potentiels}.
  Acta Sci. Math. Szeged, 1938.

\bibitem{bib:Sa}
  K.~Sato, 
  \emph{L\'evy processes and infinitely divisible distributions}.
  Cambridge Univ. Press, Cambridge, 1999.

\bibitem{bib:sw}
  R. Song, J.-M. Wu,
  \emph{Boundary Harnack principle for symmetric stable processes}.
  J. Funct. Anal. 168 (1999), 403--427

\bibitem{bib:wh}
  J.-M. Wu,
  \emph{Harmonic measures for symmetric stable processes}.  
  Studia Math.  149  (2002),  no. 3, 281--293.
 
\bibitem{bib:w}
  J.-M. Wu,
  \emph{Symmetric stable processes stay in thick sets}.
  Ann. Probab. 32 (2004), 315--336.
\end{thebibliography}
\end{document}